\documentclass{elsart}
\usepackage{amsmath}
\usepackage{amssymb}
\usepackage{epsfig}
\usepackage{graphicx}
\usepackage{subfigure}
\usepackage{amsfonts}
\usepackage{enumerate}
\usepackage{indentfirst}
\usepackage{CJK,CJKnumb}
\usepackage{color}

\newtheorem{theorem}{Theorem}
\newtheorem{example}{Example}

\newtheorem{algor}{Algorithm}

\begin{document}
\begin{CJK*}{GBK}{song}
\begin{frontmatter}

\title{Curvatures at the singular points of algebraic curves and surfaces}

\author{Chong-Jun Li\corauthref{cor}},
\ead{chongjun@dlut.edu.cn}
\author{Ren-Hong Wang}
\corauth[cor]{Corresponding author.}

\address{School of Mathematical Sciences, Dalian University
of Technology,\\ Dalian 116024, China}

\date{Revised in July 2010, the first draft in Jan 2007}

\begin{abstract}
In this paper, we study the computation of curvatures at the
singular points of algebraic curves and surfaces. The idea is to
convert the problem to compute the curvatures of the corresponding
regular parametric curves and surfaces, which have intersections
with the original curves and surfaces at the singular points. Three
algorithms are presented for three cases of plane curves, space
curves and surfaces.

\end{abstract}


\end{frontmatter}

\section{Introduction}
As we know, curvatures are important geometric information of curves
and surfaces in geometric modeling. In general, there are two ways
to represent curves and surfaces: parametric form and implicit form.
For parametric curves or surfaces, the curvature formulas were well
known \cite{doCarmo1976,Farin2002}. Recently, the curvature formulas
for implicit curves and surfaces were collected and derived in
\cite{Goldman2005}, including curvature for implicit planar curves,
curvature and torsion for implicit space curves, and mean and
Gaussian curvature for implicit surfaces. However, the curvature
formulas have been obtained only based on the regular points of the
implicit curves and surfaces. The curvature at singular point of
implicit curves or surfaces has not been discussed yet. In this
paper, three algorithms are presented for the cases of algebraic
curves and surfaces. The idea is to convert the problem to compute
the curvatures of the corresponding regular parametric curves and
surfaces, which have intersections with the original curves and
surfaces at the singular points.

For given an algebraic curve $\mathcal{C}: F(x,y) = 0$, and a
particular point $P$ of $\mathcal{C}$. If the gradient $\nabla F$
vanishes at $P$, then $P$ is a singular point of $\mathcal{C}$
\cite{Walker1978}. In view of geometry, there is only one tangent to
the curve at the regular point (where $\nabla F\neq 0$), and there
are more than one tangents (counting multiplicity) to the curve at
the singular point. The tangents are corresponding to different
branches of the curve at the singular point respectively. The number
of the tangents is the multiplicity of the singular point.
Therefore, the curvature at the singular point of the algebraic
curve can be defined by the curvature of the branch of the curve
along every tangent at the singular point. The cases of algebraic
surfaces and space algebraic curves can be discussed similarly in
the present paper.

The rest of the paper is composed of three parts for the cases of
plane algebraic curves, algebraic surfaces and space algebraic
curves, respectively. In each section, the computation of the
curvature at the singular point is presented as an algorithm. The
correctness of the algorithms are verified by the consistence with
the curvatures at the regular points. The procedure of the curvature
computation by the algorithms are illustrated by some examples
according to different cases.

\section{Curvatures at the singular points of plane algebraic curves}
Firstly, we recall the method presented in \cite{Walker1978} for
computing the tangents of an algebraic curve at its singular points.

Let $\mathcal{C}$ be a curve in the real plane $\mathbb{R}^2$,
defined by the polynomial equation $F(x,y) = 0$, and $P$ a
particular point of $\mathcal{C}$ with coordinates $(a, b)$. Let
$\mathcal{L}$ be a straight line through the point $P$, then the
parametric equations of $\mathcal{L}$ are
\begin{equation} \label{eq:PC1}
\left\{
\begin{aligned}
x&=a+\lambda t, \\
y&=b+\mu t,
\end{aligned} \right.
\end{equation}
where $\mathcal{L}$ is determined by the ratio $\lambda: \mu$. The
intersections of $\mathcal{L}$ and $\mathcal{C}$ are corresponding
to the roots of the equation $F(a+\lambda t, b+\mu t)=0$. Expanding
the left side in a Taylor series in $t$ and noticing $F(a,b)=0$, we
obtain
$$(F_x \lambda+F_y \mu)t+\frac{1}{2!}(F_{xx}\lambda^2+2F_{xy}\lambda\mu+F_{yy}\mu^2)t^2+\cdots=0,$$
where $F_x,F_y,\ldots$ are the derivatives of $F$ at $P$.

Case 1. Regular point. If not all of $F_x$, $F_y$ vanish at $P$,
then every line through $P$ has a single intersection with
$\mathcal{C}$ at $P$, with the one exception corresponding to the
value of $\lambda: \mu$ which makes $F_x \lambda+F_y \mu=0$. This
line is called the tangent to $\mathcal{C}$ at $P$, and determined
by
$$F_x(a, b)(x-a) + F_y (a, b)(y-b) = 0.$$

Case 2. Singular point. If $F_x(a,b)=F_y(a,b)=0$. Suppose that $r$
is the minimal natural number such that not all derivatives of order
$r$ of $F (x, y)$ vanish at $P$. Then every line through $P$ has at
least $r$ intersections with $\mathcal{C}$ at $P$, and precisely $r$
such lines, properly counted, have more than $r$ intersections.
These exceptional lines $\mu(x-a)-\lambda(y-b) = 0$ are called
tangents to $\mathcal{C}$ at $P$, correspond to the ratios $\lambda:
\mu$ satisfying the following equality
\begin{equation}\label{eq:tan-line}
F_{x^r}\lambda^r+\binom{r}{1}F_{x^{r-1}y}\lambda^{r-1}\mu+\cdots+{r
\choose r}F_{y^r}\mu^r=0. \end{equation} $P$ is said to be a point
of $\mathcal{C}$ of multiplicity $r$, or an $r$-fold point. A point
of $\mathcal{C}$ of multiplicity one is called a simple point or
regular point of $\mathcal{C}$. A point of multiplicity two or more
is said to be singular.

Some algebraic curves were presented in \cite{Walker1978} with
singular point at the origin as follows.

\begin{example}
$F(x,y)=x^3-x^2+y^2=0$ (see Fig.~1).
\end{example}
\begin{example}
$F(x,y)=x^3+x^2+y^2=0$ (see Fig.~2).
\end{example}
\begin{example}
$F(x,y)=x^3-y^2=0$ (see Fig.~3).
\end{example}
\begin{example}
$F(x,y)=2x^4-3x^2y+y^2-2y^3+y^4=0$ (see Fig.~4).
\end{example}
\begin{example}
$F(x,y)=x^4+x^2y^2-2x^2y-xy^2+y^2=0$ (see Fig.~5).
\end{example}
\begin{example}
$F(x,y)=(x^2+y^2)^2+3x^2y-y^3=0$ (see Fig.~6).
\end{example}
\begin{example}
$F(x,y)=(x^2+y^2)^3-4x^2y^2=0$ (see Fig.~7).
\end{example}
\begin{example}
$F(x,y)=x^6-x^2y^3-y^5=0$ (see Fig.~8).
\end{example}

\begin{figure}[!h]
\centering
\begin{minipage}[c]{0.3\textwidth}
\includegraphics[width=\textwidth]{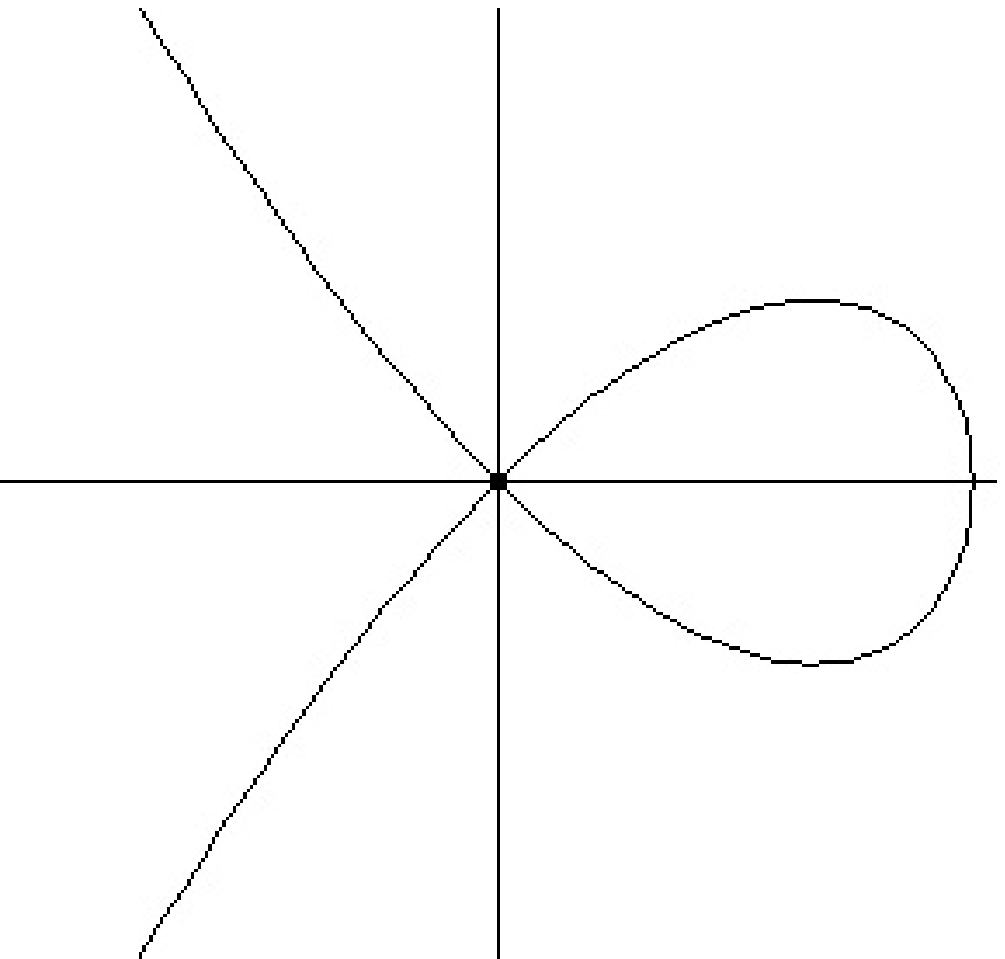}
\caption{$x^3-x^2+y^2=0$.}
\end{minipage}\hspace{3ex}
\begin{minipage}[c]{0.3\textwidth}
\includegraphics[width=\textwidth]{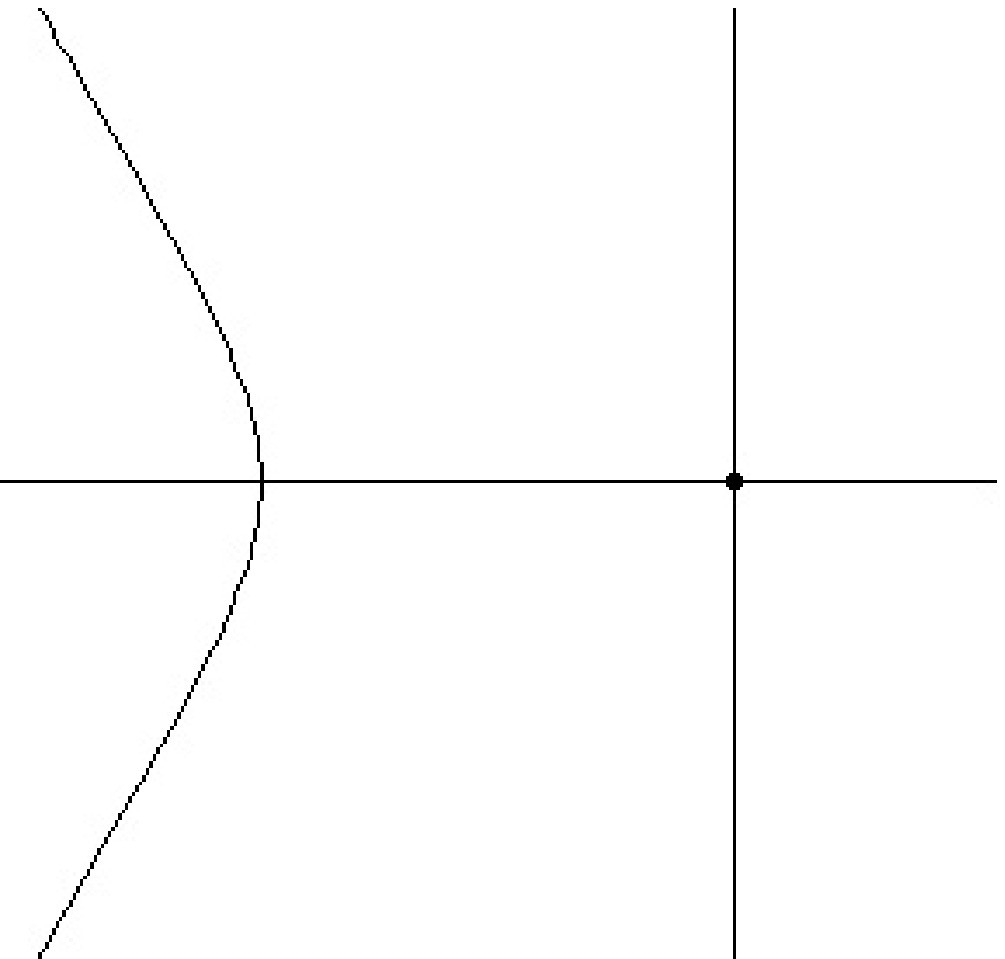}
\caption{$x^3+x^2+y^2=0$.}
\end{minipage}\hspace{3ex}
\begin{minipage}[c]{0.3\textwidth}
\includegraphics[width=\textwidth]{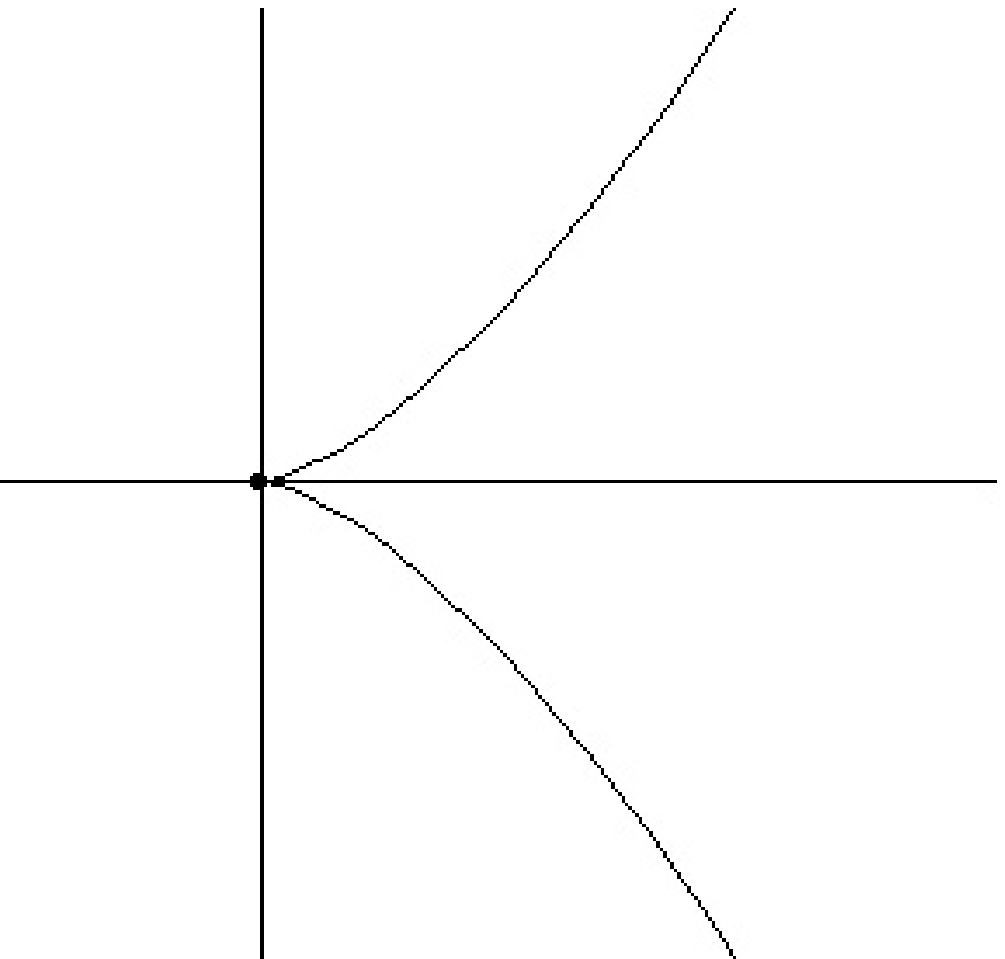}
\caption{$x^3-y^2=0$.}
\end{minipage}
\end{figure}
\begin{figure}[!h]
\centering
\begin{minipage}[t]{0.3\textwidth}
\includegraphics[width=\textwidth]{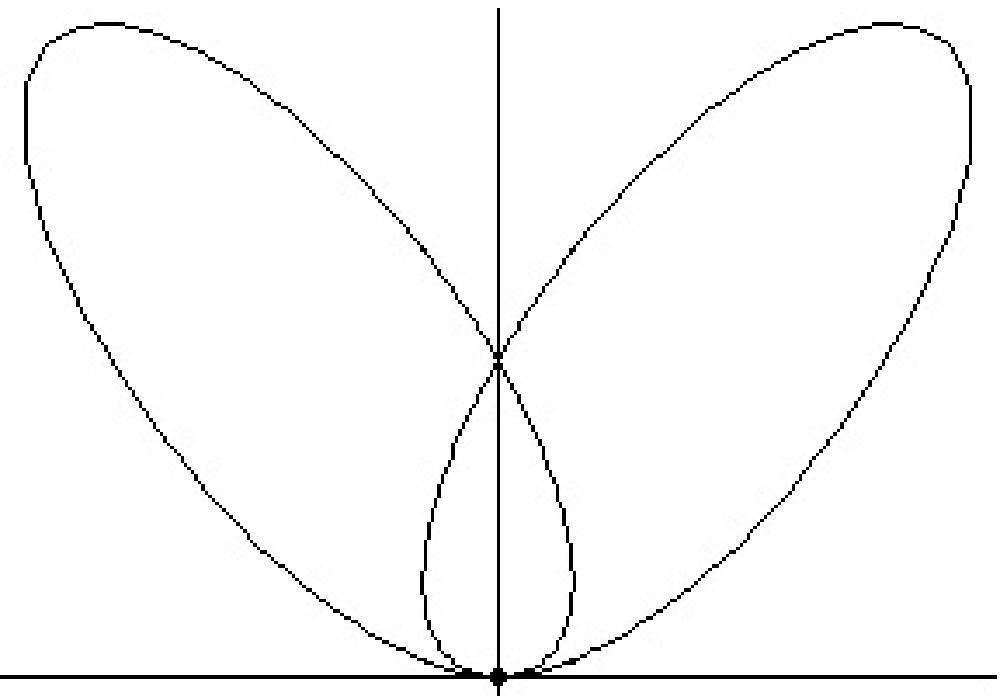}
\caption{$2x^4-3x^2y+y^2-2y^3+y^4=0$.}
\end{minipage}\hspace{3ex}
\begin{minipage}[t]{0.25\textwidth}
\includegraphics[width=\textwidth]{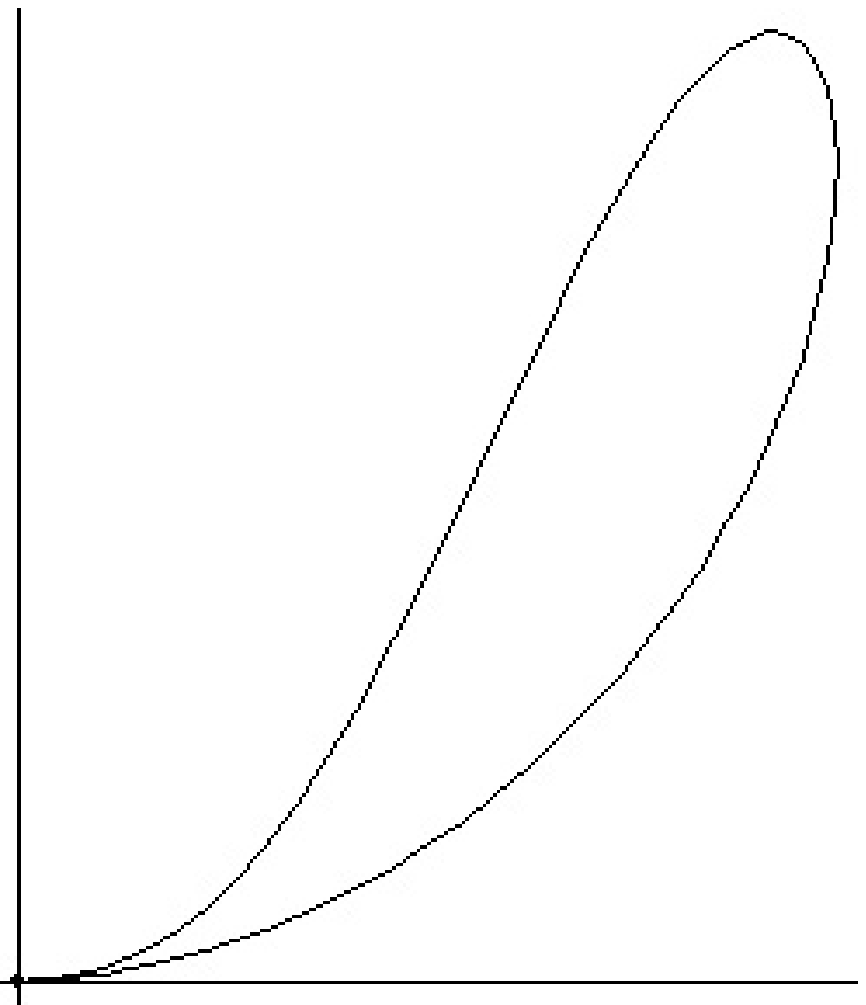}
\caption{$x^4+x^2y^2-2x^2y-xy^2+y^2=0$.}
\end{minipage}\hspace{6ex}
\begin{minipage}[t]{0.3\textwidth}
\includegraphics[width=\textwidth]{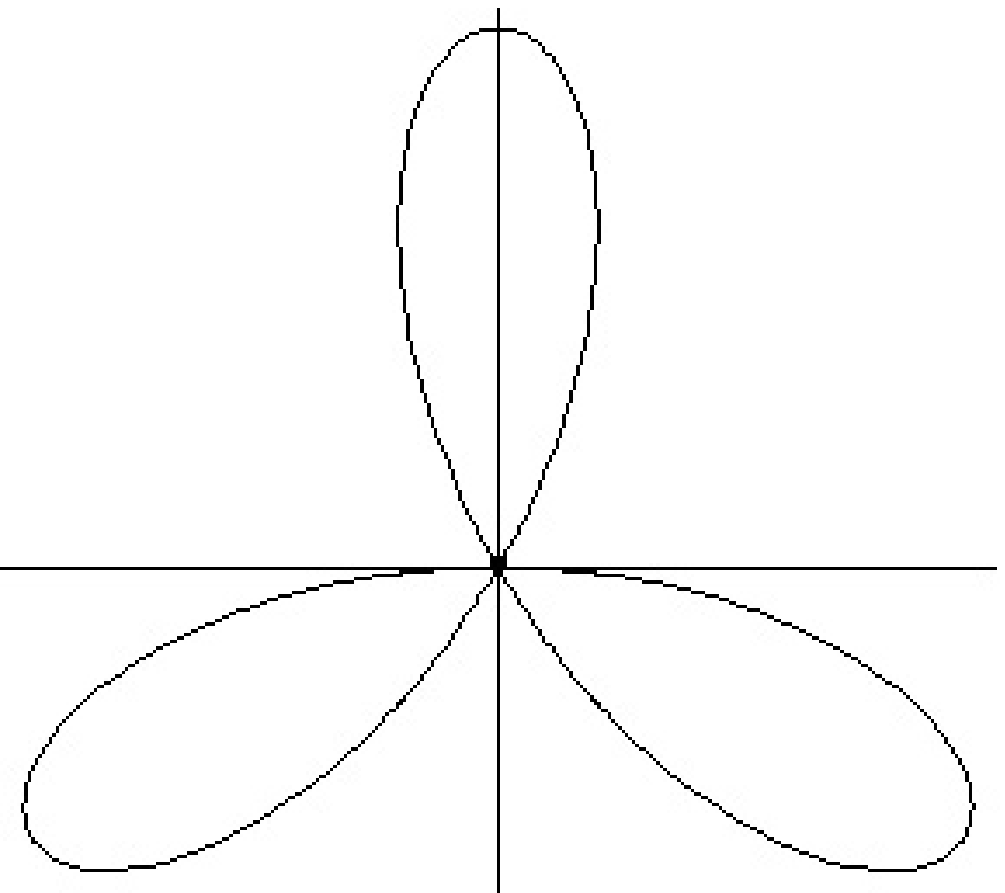}
\caption{$(x^2+y^2)^2+3x^2y-y^3=0$.}
\end{minipage}
\end{figure}
\begin{figure}[!h]
\centering
\begin{minipage}[c]{0.3\textwidth}
\includegraphics[width=\textwidth]{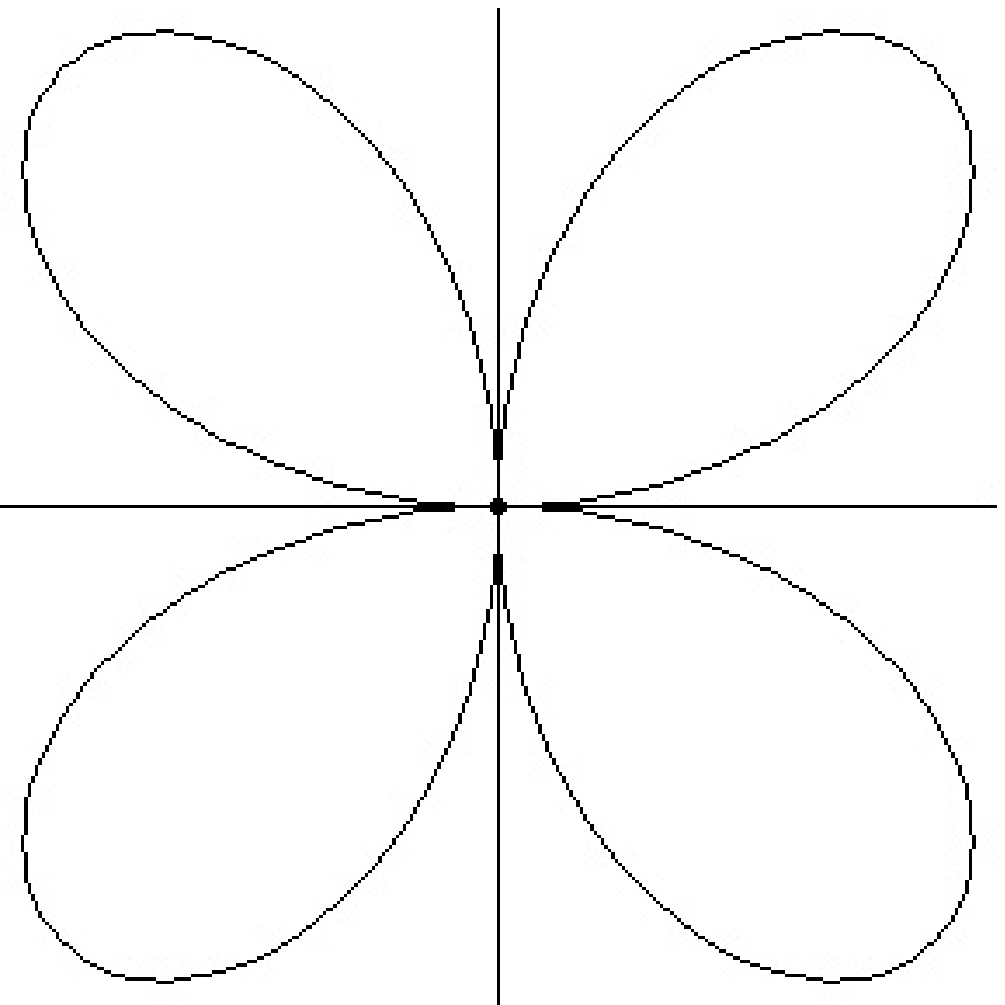}
\caption{$(x^2+y^2)^3-4x^2y^2=0$.}
\end{minipage}\hspace{3ex}
\begin{minipage}[c]{0.3\textwidth}
\includegraphics[width=\textwidth]{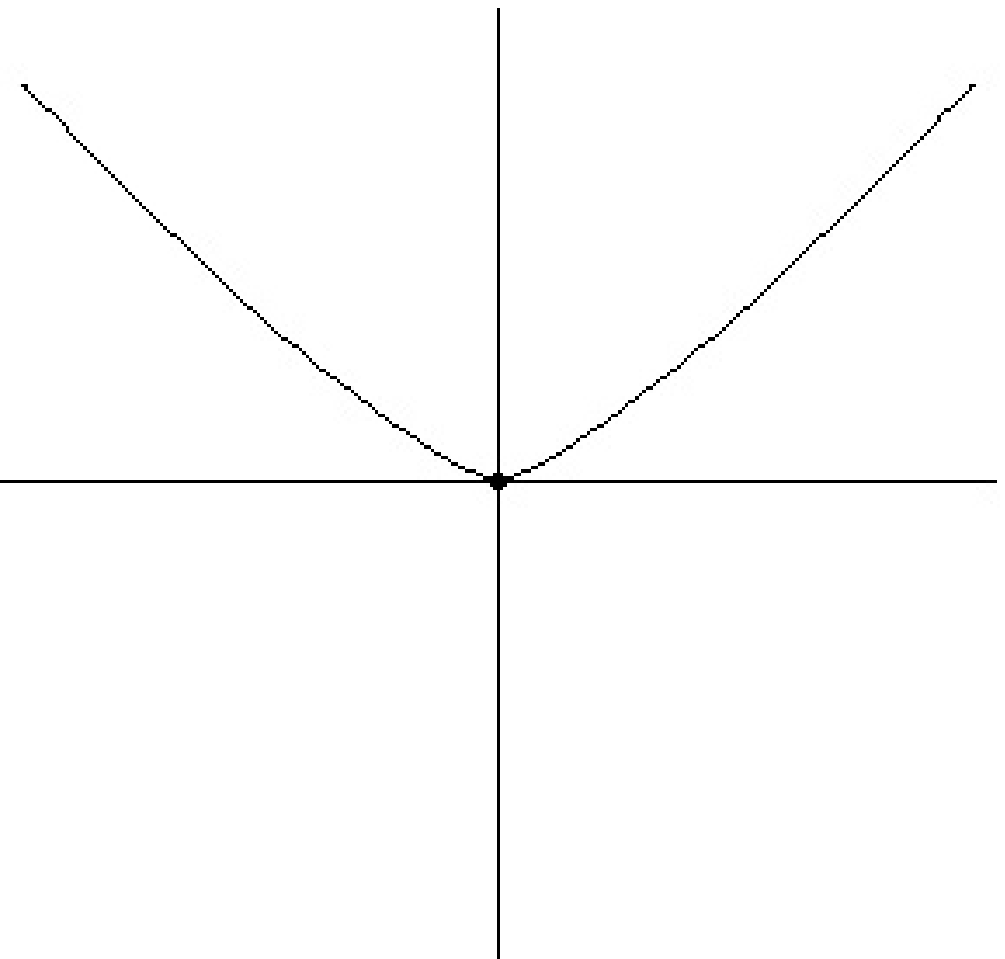}
\caption{$x^6-x^2y^3-y^5=0$.}
\end{minipage}
\end{figure}

As shown in each figure, there are different branches corresponding
with every tangent to the curve at the origin. Hence, the curvature
at the singular point can be defined as the curvatures of all
branches of the curve at the singular point. In other words, the
curvature at the singular point of the curve can be distinguished
along each tangent. By the following method, the curvatures can be
computed simultaneously.

Let $\Gamma$ be a quadratic parametric plane curve through the point
$P$ of $\mathcal{C}$, determined by
\begin{equation} \label{eq:PC2}
\mathbf{r}(t)=\left\{
\begin{aligned}
x&=a+a_1 t+\frac{1}{2}a_2 t^2, \\
y&=b+b_1 t+\frac{1}{2}b_2 t^2.
\end{aligned} \right.
\end{equation}
Then $\mathbf{r}(0)=(a,b)$, $\mathbf{r}'(0)=(a_1,b_1)$,
$\mathbf{r}''(0)=(a_2,b_2)$, and the curvature at $P=(a,b)$ of
$\Gamma$ is (\cite{doCarmo1976})
\begin{equation}\label{eq:ql1}
k=\frac{|\det(\mathbf{r}'\
\mathbf{r}'')|}{|\mathbf{r}'|^3}=\frac{|x'y''-x''y'|}{(x'^2+y'^2)^{3/2}}=\frac{|a_1b_2-a_2b_1|}{(a_1^2+b_1^2)^{3/2}}.
\end{equation}

The intersections of $\Gamma$ and $\mathcal{C}$ are corresponding to
the roots of the equation $F(x(t),y(t))=0$. Expand it by the power
of $t$ and notice that $F(a,b)=0$, then we have
\begin{eqnarray}\label{eq:Ft1}
&&F(x,y)\cr
&=&F_x(x-a)+F_y(y-b)+\frac{1}{2!}(F_{xx}(x-a)^2+2F_{xy}(x-a)(y-b)+F_{yy}(y-b)^2)+\cdots\cr
&&+\frac{1}{r!}(F_{x^r}(x-a)^r+{r\choose
1}F_{x^{r-1}y}(x-a)^{r-1}(y-b)+\cdots+{r\choose r}F_{y^r}(y-b)^r)
+\cdots\cr &=&C_1t+C_2 t^2+\cdots+C_{r}t^{r}+\cdots,
\end{eqnarray}
where $C_i$ denotes the coefficients of $t^i\ (i=1,2,\ldots)$. There
are two cases according to the singularity of $P$.
\begin{enumerate}[]
\item Case 1. If $P$ is a regular point of $\mathcal{C}$, $(F_x,F_y)\neq
0$, then
\begin{equation}\label{eq:C1}
\begin{array}{l}
C_1=F_xa_1+F_yb_1,\\
C_2=\frac{1}{2}(F_xa_2+F_yb_2)+\frac{1}{2}(F_{xx}a_1^2+2F_{xy}a_1b_1+F_{yy}b_1^2).
\end{array}
\end{equation}

\item Case 2. If $P$ is a singular point of $\mathcal{C}$ with multiplicity
$r$, then
\begin{equation}
\begin{array}{l}
C_1=0,\cr \quad\vdots \cr C_{r-1}=0,\cr
C_r=\frac{1}{r!}(F_{x^r}a_1^r+\binom{r}{1}F_{x^{r-1}y}a_1^{r-1}b_1+\cdots+{r
\choose r}F_{y^r}b_1^r).
\end{array}
\end{equation}
\end{enumerate}

By $C_r=0$, we obtain the tangent vectors $\mathbf{r}'(0)=(a_1,b_1)$
same to Eq.~(\ref{eq:tan-line}). It means that the parametric curve
$\Gamma$ has all the same tangents with the algebraic curve
$\mathcal{C}$ at the singular point $P$, and they have the
intersection point $P$ with multiplicity more than $r$. Note that
the root $\mathbf{r}'(0)=(a_1,b_1)\neq 0$, i.e., $P$ is a regular
point of the parametric curve $\Gamma$. By Frenet equations, a
parametric curve is determined by the tangent, curvature and torsion
at its regular point. Hence, when $\Gamma$ and $\mathcal{C}$ have
the intersection point $P$ with multiplicity much more than $r$,
they will have same curvatures at $P$ along all tangents.

Consequently, the computation of the curvature at the singular point
of the algebraic curve $\mathcal{C}$ along the tangent $(a_1,b_1)$
is converted to compute the curvature at the regular point of the
quadratic parametric curve $\Gamma$ along the same tangent. We
conclude the procedure as the following algorithm.

\begin{algor}[Curvatures at the singular point of a plane algebraic
curve]\label{Alg:ql1}\hspace{1cm}\\ For given a plane algebraic
curve $\mathcal{C}: F(x,y)=0$, and $P=(a,b)$ is an $r$-fold point of
$\mathcal{C}$. Let $\Gamma$ be a parametric quadratic curve through
the point $P$ defined by Eq.~(\ref{eq:PC2}), then the curvature at
$P$ along every tangent to $\mathcal{C}$ is determined by the
quadratic curve $\Gamma$, which has intersection with $\mathcal{C}$
at $P$ with multiplicity of $l>r+1$.
\begin{enumerate}[1)]
\item
The $r$ tangents $(a_1,b_1)$ to $\mathcal{C}$ at $P$ are determined
by $C_r=0$, where $C_r$ satisfying Eq.~(\ref{eq:Ft1}).

\item Substitute the above $a_1,\ b_1$ into Eq.~(\ref{eq:Ft1}), if $l$ is the minimal natural number such
that $C_r=\cdots=C_{l-1}=0,\ C_l\neq0$, then the curvature of
$\mathcal{C}$ at $P$ along the tangent $(a_1,b_1)$ equals to the
curvature of $\Gamma$ at $P$, determined by $C_l=0$. The curvature
is computed by formula~(\ref{eq:ql1}).
\end{enumerate}
\end{algor}

In order to verify the correctness of Algorithm \ref{Alg:ql1}, we
firstly prove that the curvature obtained by Algorithm~\ref{Alg:ql1}
is same to that by formula in \cite{Goldman2005} if $P$ is regular.

\begin{theorem}
If $P$ is a regular point of $\mathcal{C}$, the curvature at $P$ by
Algorithm \ref{Alg:ql1} is equivalent to the curvature formula for
regular implicit plane curves,
\begin{equation}
k=\frac{|{\rm Tan}(F)*H(F)*{\rm Tan}(F)^T|}{|\nabla
F|^3}=\frac{|(-F_y\ F_x)\binom{F_{xx}\ F_{xy}}{F_{yx} \
F_{yy}}\binom{-F_y}{F_x}|}{(F_x^2+F_y^2)^{3/2}}.
\end{equation}
\end{theorem}
\begin{pf*}{Proof}
The multiplicity of $P$ is $r=1$. By $C_1=0$ in Eq.~(\ref{eq:C1}),
then $a_1: b_1=-F_y: F_y$. Without loss of generality, we may assume
$a_1=-F_y,\ b_1=F_x$, then
$$C_2=\frac{1}{2}(b_1a_2-a_1b_2)+\frac{1}{2}(F_{xx}F_y^2-2F_{xy}F_xF_y+F_{yy}F_y^2).$$
By $C_2=0$ and curvature formula (\ref{eq:ql1}) for parametric plane
curve $\Gamma$, we obtain the curvature at $P$ of $\mathcal{C}$ as
$$k=\frac{|a_1b_2-a_2b_1|}{(a_1^2+b_1^2)^{3/2}}=\frac{|F_{xx}F_y^2-2F_{xy}F_xF_y+F_{yy}F_y^2|}{(F_x^2+F_y^2)^{3/2}}
=\frac{|(-F_y\ F_x)\binom{F_{xx}\ F_{xy}}{F_{yx} \
F_{yy}}\binom{-F_y}{F_x}|}{(F_x^2+F_y^2)^{3/2}}.$$ The proof is
completed.
\end{pf*}

Next, we take an example of a reducible curve to illustrate the
computation of curvatures at the singular point by Algorithm
\ref{Alg:ql1}.
\begin{example}{\rm  $F(x,y)=(x-y)(x^2+y^2-2Rx)$.
There are two singular points $(0,0)$ and $(R,R)$ of the curve
$F(x,y)=0$. We only focus on the case of $P=(0,0)$, the other one
can be discussed similarly. As shown in Fig.~\ref{fig:ex9}, it is
clear that the curve is composed of a straight line $x=y$ and a
circle with a radius of $R$ crossing at the origin. Hence, the two
tangent vectors are $(1,1)$ and $(0,1)$, and the two curvatures are
0 and $1/R$ corresponding to the straight line and the circle,
respectively.
\begin{figure}[!h]
\centering
\includegraphics[width=0.3\textwidth]{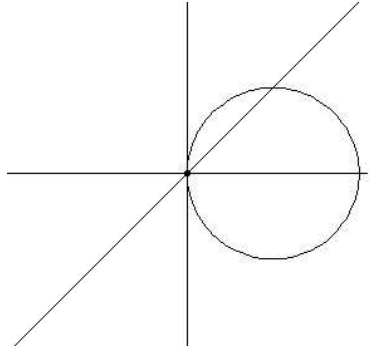}
\caption{$(x-y)(x^2+y^2-2Rx)=0$.}\label{fig:ex9}
\end{figure}

By Algorithm \ref{Alg:ql1}, let
$$\mathbf{r}(t)=\left\{
\begin{aligned}
x(t)&=a_1 t+\frac{1}{2}a_2 t^2, \\
y(t)&=b_1 t+\frac{1}{2}b_2 t^2,
\end{aligned} \right.$$
then
$$F(x,y)=-2a_1(a_1-b_1)Rt^2+((a_1-b_1)(a_1^2+b_1^2-a_2R)-a_1(a_2-b_2)R)t^3+\cdots.$$
By Eq.~(\ref{eq:Ft1}), $C_2=-2a_1(a_1-b_1)R\neq 0$, hence, $P=(0,0)$
is a 2-fold point (or double point). Let $C_2=0$, then we have two
solutions $a_1=b_1$ and $a_1=0$ corresponding to the two following
cases respectively.

\begin{enumerate}[1)]
\item When $a_1=b_1\neq 0$, then $C_2=0$, $C_3=-a_1(a_2-b_2)$, and the tangent is $(1,1)$.
Let $C_3=0$ because the order of the singularity is two, then
$a_2=b_2$. By Eq.~(\ref{eq:ql1}), we obtain the curvature $k=0$
along the tangent $(1,1)$.

\item When $a_1=0,\ b_1\neq 0$, then $C_2=0$, $C_3=-b_1(b_1^2-a_2R)$, and the tangent is $(0,1)$.
Let $C_3=0$, then $a_2=b_1^2/R$, the curvature $k=|a_2/b_1^2|=1/R$
along the tangent $(0,1)$.
\end{enumerate}
}\end{example}

Further, we compute the curvatures at the origin of the curves in
Examples 1-8 by Algorithm \ref{Alg:ql1}, as shown in Figures 1-8.
\begin{enumerate}[Ex.~1]
\item $F(x,y)=x^3-x^2+y^2=0$. $C_1=0$, $C_2=-a_1^2+b_1^2$,
$C_3=a_1^3-a_1a_2+b_1b_2$, and the multiplicity $r=2$.
\begin{enumerate}[1)]
\item When $a_1=b_1\neq 0$, then $C_2=0$, $C_3=-b_1(a_2-b_1^2-b_2)$.
Let $a_2=b_1^2+b_2$, then $C_3=0$, $k=\sqrt{2}/4$ along the tangent
$(1,1)$.

\item When $a_1=-b_1\neq 0$, then $C_2=0$, $C_3=b_1(a_2-b_1^2+b_2)$.
Let $a_2=b_1^2-b_2$, then $C_3=0$, $k=\sqrt{2}/4$ along the tangent
$(1,-1)$.
\end{enumerate}

\item $F(x,y)=x^3+x^2+y^2=0$. $C_1=0$, $C_2=a_1^2+b_1^2$,
$C_3=a_1^3+a_1a_2+b_1b_2$, and the multiplicity $r=2$.
\begin{enumerate}[1)]
\item When $a_1=1,\ b_1=i$, then $C_2=0$, $C_3=1+a_2+b_2 i$.
Let $a_2=-1-b_2i$, then $C_3=0$, $k=|a_1b_2-a_2b_1|/0=|i|/0=\infty$
along the tangent $(1,i)$.

\item When $a_1=1,\ b_1=-i$, then $C_2=0$, $C_3=1+a_2-b_2 i$.
Let $a_2=-1+b_2i$, then $C_3=0$, $k=|a_1b_2-a_2b_1|/0=|-i|/0=\infty$
along the tangent $(1,-i)$.
\end{enumerate}

\item $F(x,y)=x^3-y^2=0$. $C_1=0$, $C_2=-b_1^2$, $C_3=a_1^3-b_1b_2$, and the multiplicity $r=2$.
Since $b_1=0$ is a double root of $C_2=0$, there is a double tangent
$(1,0)$ at the singular point.

When $b_1=0$, then $C_2=0$, $C_3=a_1^3$. We want $C_3$ to be zero
because the order of singularity is two, then we have $a_1=0$.
However, it means that the origin is also a singular point of the
parametric curve $\mathbf{r}(t)$ as
$\mathbf{r}'(t)=(a_1,b_1)=(0,0)$. In fact, the parametric curve is
$$\mathbf{r}(t)=\left\{
\begin{aligned}
x(t)&=\frac{1}{2}a_2 t^2, \\
y(t)&=\frac{1}{2}b_2 t^2,
\end{aligned} \right.$$
i.e., it is a ray with the end point $P=(0,0)$. Because the meaning
of curvature is to measure the speed of rotation of the tangent, we
can define the curvature of the ray at the end point by infinity.

On another hand, by solving the equation $F(x,y)=x^3-y^2=0$, we
obtain the two components of the curve, $y=\sqrt{x^3}$ and
$y=-\sqrt{x^3}$. The condition $x\geq 0$ is necessary such that the
curve is real. Hence, the singular point $P$ is the end point of the
algebraic curve. By Eq.~(\ref{eq:ql1}), the curvature of the
component $y=\sqrt{x^3}$ is $$k(x)=\frac{6}{\sqrt{x}(4+9x)^{3/2}}.$$
When $x$ tends to zero, the limit of $k(x)$ is infinity. Hence, the
curvature of the curve $F(x,y)=x^3-y^2=0$ at the origin is $\infty$
along the double tangent (1,0).


\item $F(x,y)=2x^4-3x^2y+y^2-2y^3+y^4=0$. $C_1=0$, $C_2=b_1^2$, $C_3=-b_1(3a_1^2+2b_1^2-b_2)$, and the multiplicity $r=2$.
When $b_1=0$, then $C_2=C_3=0$,
$C_4=\frac{1}{4}(2a_1^2-b_2)(4a_1^2-b_2)$, and there is a double
tangent $(1,0)$ at the singular point.
\begin{enumerate}[1)]
\item Let $b_2=2a_1^2$, then $C_4=0$, $k=|b_2|/a_1^2=2$.
\item Let $b_2=4a_1^2$, then $C_4=0$, $k=|b_2|/a_1^2=4$.
\end{enumerate}

\item $F(x,y)=x^4+x^2y^2-2x^2y-xy^2+y^2=0$. $C_1=0$, $C_2=b_1^2$, $C_3=-b_1(2a_1^2+a_1b_1^2-b_2)$, and the multiplicity $r=2$.
When $b_1=0$, then $C_2=C_3=0$, $C_4=\frac{1}{4}(2a_1^2-b_2)^2$. Let
$b_2=2a_1^2$, then $C_4=0$, $k=|b_2|/a_1^2=2$ along the double
tangent $(1,0)$.

\item $F(x,y)=(x^2+y^2)^2+3x^2y-y^3=0$. $C_1=C_2=0$, $C_3=b_1(3a_1^2-b_1^2)$, and the multiplicity $r=3$.
\begin{enumerate}[1)]
\item When $b_1=0$, then $C_3=0$, $C_4=\frac{1}{2}a_1^2(2a_1^2+3b_2)$.
Let $b_2=-\frac{2}{3}a_1^2$, then $C_4=0$,
$k=|b_2|/a1^2=\frac{2}{3}$ along the tangent $(1,0)$.

\item When $a_1=1,\ b_1=\sqrt{3}$, then $C_3=0$, $C_4=16+3\sqrt{3}a_2-3b_2$.
Let $b_2=\frac{16}{3}+\sqrt{3}a_2$, then $C_4=0$, $k=\frac{2}{3}$
along the tangent $(1,\sqrt{3})$.

\item When $a_1=1,\ b_1=-\sqrt{3}$, then $C_3=0$, $C_4=16-3\sqrt{3}a_2-3b_2$.
Let $b_2=\frac{16}{3}-\sqrt{3}a_2$, then $C_4=0$, $k=\frac{2}{3}$
along the tangent $(1,-\sqrt{3})$.
\end{enumerate}

\item $F(x,y)=(x^2+y^2)^3-4x^2y^2=0$. $C_1=C_2=C_3=0$,
$C_4=-4a_1^2b_1^2$, and the multiplicity $r=4$.
\begin{enumerate}[1)]
\item When $a_1=0$, then $C_4=C_5=0$, $C_6=-b_1^2(a_2-b_1^2)(a_2+b_1^2)$.
Let $a_2=b_1^2$ or $a_2=-b_1^2$, then $C_6=0$, $k=|a_2|/b_1^2=1$
along the double tangent $(0,1)$.

\item When $b_1=0$, then $C_4=C_5=0$, $C_6=a_1^2(a_1^2-b_2)(a_1^2+b_2)$.
Let $b_2=a_1^2$ or $b_2=-a_1^2$, then $C_6=0$, $k=|b_2|/a_1^2=1$
along the double tangent $(1,0)$.
\end{enumerate}

\item $F(x,y)=x^6-x^2y^3-y^5=0$. $C_1=C_2=C_3=C_4=0$,
$C_5=-b_1^3(a_1^2+b_1^2)$, and the multiplicity $r=5$.
\begin{enumerate}[1)]
\item When $b_1=0$, then $C_5=0$, $C_6=a_1^6$. It is same to the
case in Ex. 3. Then the curvature of the curve $F(x,y)=0$ at the
origin is $\infty$ along the tangent (1,0) of multiplicity 3.


\item When $a_1=1,\ b_1=i$, then $C_5=0$, $C_6=1+a_2i-b_2$.
Let $b_2=1+a_2i$, then $C_6=0$, $k=|a_1b_2-a_2b_1|/0=1/0=\infty$
along the tangent $(1,i)$.

\item When $a_1=1,\ b_1=-i$, then $C_5=0$, $C_6=1-a_2i-b_2$.
Let $b_2=1-a_2i$, then $C_6=0$, $k=|a_1b_2-a_2b_1|/0=1/0=\infty$
along the tangent $(1,-i)$.
\end{enumerate}
\end{enumerate}


The tangents and curvatures at the origin of the curves in Examples
1-8 are concluded in the following table.
\begin{table}[!h]
\centering \caption{Curvatures at the origin on the curves in
Examples 1-8.} \label{lab:1}
\begin{tabular}{|c||c|c||c|c||c|c||c|c||c|c|}
\hline Ex. & Tan.1 & Cur.1 & Tan.2 & Cur.2 & Tan.3 & Cur.3 & Tan.4 & Cur.4&Tan.5 & Cur.5\\
\hline 1& $(1,1)$ & $\sqrt{2}/4$ & $(1,-1)$ & $\sqrt{2}/4$ &&&&&&\\
\hline 2& $(1,i)$ & $\infty$ & $(1,-i)$ & $\infty$ &&&&&&\\
\hline 3& $(1,0)$ & $\infty$ & $(1,0)$ & $\infty$ &&&&&&\\
\hline 4& $(1,0)$ & $2$ & $(1,0)$ & $4$ &&&&&&\\
\hline 5& $(1,0)$ & $2$ & $(1,0)$ & $2$ &&&&&&\\
\hline 6& $(1,0)$ & $2/3$ & $(1,\sqrt{3})$ & $2/3$ & $(1,-\sqrt{3})$ & $2/3$ &&&&\\
\hline 7& $(1,0)$ & $1$ & $(1,0)$ & $1$ & $(0,1)$ & $1$ & $(0,1)$ & $1$&&\\
\hline 8& $(1,0)$ & $\infty$ & $(1,0)$ & $\infty$ & $(1,0)$ & $\infty$ & $(1,i)$ & $\infty$ & $(1,-i)$ & $\infty$\\
\hline
\end{tabular}
\end{table}

\section{Curvatures at the singular points of algebraic surfaces}
For given an algebraic surface $\mathcal{S}:\ F(x,y,z)=0$, we
discuss the Gauss and mean curvatures at the singular point $P$,
where gradient $\nabla F=0$. Let $\Sigma$ be a quadratic parametric
surface through the point $P=(a,b,c)$ on $\mathcal{S}$, denoted by
\begin{equation} \label{eq:PS}
\mathbf{r}(s,t)=\left\{
\begin{aligned}
x(s,t)&=a+a_1 s+a_2 t+\frac{1}{2}(a_3 s^2+ 2a_4st+ a_5t^2), \\
y(s,t)&=b+b_1 s+b_2 t+\frac{1}{2}(b_3 s^2+ 2b_4st+ b_5t^2), \\
z(s,t)&=c+c_1 s+c_2 t+\frac{1}{2}(c_3 s^2+ 2c_4st+ c_5t^2).
\end{aligned} \right.
\end{equation}
Then $$\mathbf{r}(0,0)=(a,b,c),\ \mathbf{r}_s(0,0)=(a_1,b_1,c_1),\
\mathbf{r}_t(0,0)=(a_2,b_2,c_2),$$
$$\mathbf{r}_{ss}(0,0)=(a_3,b_3,c_3),\
\mathbf{r}_{st}(0,0)=(a_4,b_4,c_4),\
\mathbf{r}_{tt}(0,0)=(a_5,b_5,c_5).$$ The normal vector to the
surface is perpendicular to the tangent vectors $\mathbf{r}_s$ and
$\mathbf{r}_t$. Therefore, the unit normal is given by
$\mathbf{n}=\mathbf{n}(s,t)=\frac{\mathbf{r}_s\times
\mathbf{r}_t}{|\mathbf{r}_s\times \mathbf{r}_t|}$. The first and
second fundamental forms of the surface are given by the following
matrices
$$I=\binom{E\ F}{F\ G}=\binom{\mathbf{r}_s\cdot \mathbf{r}_s \quad \mathbf{r}_s\cdot \mathbf{r}_t}
{\mathbf{r}_t\cdot \mathbf{r}_s \quad \mathbf{r}_t\cdot
\mathbf{r}_t},$$
$$II=\binom{L\ M}{M\ N}=\binom{\mathbf{r}_{ss}\cdot \mathbf{n} \quad \mathbf{r}_{st}\cdot \mathbf{n}}
{\mathbf{r}_{ts}\cdot \mathbf{n} \quad \mathbf{r}_{tt}\cdot
\mathbf{n}}.$$ Then the Gauss and mean curvatures are computed by
\begin{eqnarray}
&&K_G=\frac{LN-M^2}{EG-F^2}, \label{eq:qlg}\\
&&K_M=\frac{EN-2FM-GL}{2(EG-F^2)}.\label{eq:qlm}
\end{eqnarray}

The intersections of $\Sigma$ and $\mathcal{S}$ are corresponding to
the roots of the equation $F(x(s,t),y(s,t)$, $z(s,t))=0$. Expand it
by the power of $s$ and $t$ and notice that $F(a,b,c)=0$, then we
have
\begin{eqnarray}\label{eq:Fst}
&&F(x,y,z)\cr
&=&F_x(x-a)+F_y(y-b)+F_z(z-c)+\frac{1}{2!}(F_{xx}(x-a)^2+F_{yy}(y-b)^2+F_{zz}(z-c)^2\cr
&&+ 2F_{xy}(x-a)(y-b)+2F_{yz}(y-b)(z-c)+2F_{xz}(x-a)(z-c))+\cdots\cr
&&+
\sum_{i+j+k=r}\frac{r!}{i!j!k!}F_{x^iy^jz^k}(x-a)^i(y-b)^j(z-c)^k+\cdots\cr
&=&C_{1,0}s+C_{0,1} t+C_{2,0} s^2+C_{1,1} st +C_{0,2} t^2+\cdots,
\end{eqnarray}
where $C_{i,j}\ (i,j=1,2,\ldots)$ denote the coefficients of
$s^it^j$. Denote the vector of coefficients of degree $r$ by
$C_r=(C_{r,0},C_{r-1,1},\ldots,C_{0,r})$, $r=1,2,\ldots$. If
$C_1=\cdots=C_{r-1}=0$ and $C_r\neq0$, then $(0,0)$ is a root of
$F(s,t)$ of multiplicity $r$.

By Eq.~(\ref{eq:PS}), the tangent plane at $P$ to the parametric
surface $\Sigma$ is determined by its normal vector
$\mathbf{r}_s\times \mathbf{r}_t=(a_1,b_1,c_1)\times(a_2,b_2,c_2)$.
If $r$ is the minimal natural number such that not all partial
derivatives of order $r$ of $F(x, y)$ vanish at $P$. Then each
surface $\Sigma$ through $P$ has at least $r$ intersections with
$\mathcal{S}$ at $P$. Let $C_r=0$, then $\Sigma$ and $\mathcal{S}$
will have more than $r$ intersections. Notice that all the
coefficients $C_{r,0},C_{r-1,1},\ldots,C_{0,r}$ of degree $r$ are
only depending on $a_1,b_1,c_1,a_2,b_2,c_2$, hence there exist $r$
common tangent planes (properly computed, and let
$(a_1,b_1,c_1)\times(a_2,b_2,c_2)\neq 0$) of $\Sigma$ and
$\mathcal{S}$ at $P$, corresponding to the root of
\begin{equation}\label{eq:tan-plane}
C_r=(C_{r,0},C_{r-1,1},\ldots,C_{0,r})=0. \end{equation} $P$ is said
to be a point of $\mathcal{S}$ of multiplicity $r$. A point of
$\mathcal{S}$ of multiplicity one is called a simple point or
regular point of $\mathcal{S}$. A point of multiplicity two or more
is said to be singular.

\begin{enumerate}[1)]
\item If $P$ is a regular point of $\mathcal{S}$, $\nabla
F=(F_x,F_y,F_z)\neq 0$, then
\begin{eqnarray}\label{eq:C2}
\begin{array}{ll}
C_{1,0}&=F_xa_1+F_yb_1+F_zc_1,\\
C_{0,1}&=F_xa_2+F_yb_2+F_zc_2,\\
C_{2,0}&=\frac{1}{2}(F_xa_3+F_yb_3+F_zc_3+F_{xx}a_1^2+F_{yy}b_1^2+F_{zz}c_1^2\\&+2F_{xy}a_1b_1+2F_{yz}b_1c_1+2F_{xz}a_1c_1),\\
C_{1,1}&=F_xa_4+F_yb_4+F_zc_4+F_{xx}a_1a_2+F_{yy}b_1b_2+F_{zz}c_1c_2\\&+F_{xy}(a_1b_2+a_2b_1)+F_{yz}(b_1c_2+b_2c_1)+F_{xz}(a_1c_2+a_2c_1),\\
C_{0,2}&=\frac{1}{2}(F_xa_5+F_yb_5+F_zc_5+F_{xx}a_2^2+F_{yy}b_2^2+F_{zz}c_2^2\\&+2F_{xy}a_2b_2+2F_{yz}b_2c_2+2F_{xz}a_2c_2).\\
\end{array}
\end{eqnarray}

\item If $P$ is a singular point of $\mathcal{S}$ of multiplicity $r$,
then $C_1=\cdots= C_{r-1}=0$, $C_r\neq 0$.
\end{enumerate}

Similar to the case of plane algebraic curves, there are $r$ tangent
planes to the surface at the singular point. Then the curvature of
the surface at the singular point can be defined by the curvature of
the branch of the surface along every tangent plane or normal
vector. Note that the normal vector
$(a_1,b_1,c_1)\times(a_2,b_2,c_2)\neq 0$ which is determined by
$C_r=0$, i.e., $P$ is a regular point of the parametric surface
$\Sigma$. Hence, if $\Sigma$ and $\mathcal{S}$ have the intersection
point $P$ of multiplicity much more than $r$, they will have same
curvatures at $P$. The computing procedure can be concluded as
follows.

\begin{algor}[Curvatures at the singular point of an algebraic
surface]\label{Alg:Sql} For given an algebraic surface $\mathcal{S}:
F(x,y,z)=0$, $P=(a,b,c)$ is an $r$-fold point of $\mathcal{S}$. Let
$\Sigma$ be a parametric quadratic surface through the point $P$
defined by Eq.~(\ref{eq:PS}), then the curvature at $P$ along each
tangent plane to $\mathcal{S}$ is determined by the quadratic
surface $\Sigma$, which has intersection with $\mathcal{S}$ at $P$
with multiplicity of $l>r+1$.
\begin{enumerate}[1)]
\item
the $r$ normal vectors $(a_1,b_1,c_1)\times(a_2,b_2,c_2)$ of the
tangent planes to $\mathcal{S}$ at $P$ are determined by equation
$C_r=0$, where $C_r$ satisfying Eq.~(\ref{eq:Fst}).

\item substitute the above $a_1,b_1,c_1,a_2,b_2,c_2$ into Eq.~(\ref{eq:Fst}), if $l$ is the minimal natural number such
that $C_r=\cdots=C_{l-1}=0,\ C_l\neq0$, then the curvature of
$\mathcal{S}$ at $P$ along the normal vector
$(a_1,b_1,c_1)\times(a_2,b_2,c_2)$ equals to the curvature of
$\Sigma$ at $P$, which determined by equation $C_l=0$. The
curvatures are computed by formulas (\ref{eq:qlg}) and
(\ref{eq:qlm}).
\end{enumerate}
\end{algor}

By Algorithm \ref{Alg:Sql}, if $P$ is a regular point of
$\mathcal{S}$, we have the same results with \cite{Goldman2005} as
follows.
\begin{theorem}
If $P$ is a regular point of an algebraic surface $\mathcal{S}$, by
Algorithm \ref{Alg:Sql}, the Gauss and mean curvatures at $P$ of
$\mathcal{S}$ are
\begin{eqnarray}
K_G&=&\frac{LN-M^2}{EG-F^2}=\frac{\nabla F *H^*(F)* \nabla
F^T}{|\nabla F|^4}, \\
K_M&=&\frac{EN-2FM+GL}{2(EG-F^2)}=\frac{\nabla F*H^*(F)*\nabla
F^T-|\nabla F|^2Trace(H(F))}{2|\nabla F|^3}.
\end{eqnarray}
\end{theorem}
\begin{pf*}{Proof}
The multiplicity of $P$ is $r=1$. By $C_{1,0}=C_{0,1}=0$, we obtain
$\mathbf{r}_s\times \mathbf{r}_t=\lambda \nabla F$, for a constant
$\lambda\neq 0$. Let $\lambda=1$ for convenience.

Denote the hessian by
$$H=H(F)=\left( \begin{matrix}
F_{xx}&F_{xy}&F_{xz}\cr F_{yx}&F_{yy}&F_{yz}\cr F_{zx}&F_{zy}&F_{zz}
\end{matrix}\right)=\nabla(\nabla F),$$
and the adjoint of the hessian $H^*=H^*(F)$. Then by
$C_{2,0}=C_{1,1}=C_{0,2}=0$, we have
$$\begin{array}{ll}\mathbf{r}_{ss}\cdot \nabla F&=-\mathbf{r}_s
*H*\mathbf{r}_s^T,\cr \mathbf{r}_{st}\cdot \nabla
F&=-\mathbf{r}_s*H*\mathbf{r}_t^T=-\mathbf{r}_t *H*
\mathbf{r}_s^T,\cr \mathbf{r}_{tt}\cdot \nabla F&=-\mathbf{r}_t
*H* \mathbf{r}_t^T.
\end{array}
$$

Note that $\mathbf{n}=\frac{\nabla F}{|\nabla F|}$, then
\begin{eqnarray*}
EG-F^2&=&|\mathbf{r}_s\times\mathbf{r}_t|^2=\lambda^2|\nabla F|^2,
\end{eqnarray*}
\begin{eqnarray*}
LN-M^2&=&(\mathbf{r}_{ss}\cdot \mathbf{n})(\mathbf{r}_{tt}\cdot
\mathbf{n}) -(\mathbf{r}_{st}\cdot
\mathbf{n})(\mathbf{r}_{ts}\cdot\mathbf{n}) \cr &=&\frac{1}{|\nabla
F|^2}((\mathbf{r}_{ss}\cdot \nabla F)(\mathbf{r}_{tt}\cdot \nabla F)
-(\mathbf{r}_{st}\cdot \nabla F)(\mathbf{r}_{ts}\cdot\nabla F))\cr
&=& \frac{1}{|\nabla F|^2}((\mathbf{r}_s*
H\cdot\mathbf{r}_s)(\mathbf{r}_t* H\cdot\mathbf{r}_t)-(\mathbf{r}_s*
H\cdot \mathbf{r}_t)(\mathbf{r}_t* H\cdot \mathbf{r}_s))\cr
&=&\frac{1}{|\nabla F|^2} (\mathbf{r}_s* H\times\mathbf{r}_t*
H)\cdot(\mathbf{r}_s\times\mathbf{r}_t) \cr &=&\frac{1}{|\nabla
F|^2}((\mathbf{r}_s\times \mathbf{r}_t)*
H^*)\cdot(\mathbf{r}_s\times\mathbf{r}_t),
\end{eqnarray*}
\begin{eqnarray*}
&& EN-2FM+GL \cr
&=&(\mathbf{r}_s\cdot\mathbf{r}_s)(\mathbf{r}_{tt}\cdot\mathbf{n})
-2(\mathbf{r}_s\cdot\mathbf{r}_t)(\mathbf{r}_{st}\cdot\mathbf{n})+(\mathbf{r}_t\cdot\mathbf{r}_t)(\mathbf{r}_{ss}\cdot\mathbf{n})\cr
&=&\frac{1}{|\nabla
F|}((\mathbf{r}_s\cdot\mathbf{r}_s)(\mathbf{r}_{tt}\cdot\nabla F)
-2(\mathbf{r}_s\cdot\mathbf{r}_t)(\mathbf{r}_{st}\cdot\nabla
F)+(\mathbf{r}_t\cdot\mathbf{r}_t)(\mathbf{r}_{ss}\cdot\nabla F))
\cr &=&\frac{-1}{|\nabla
F|}((\mathbf{r}_s\cdot\mathbf{r}_s)(\mathbf{r}_{t}*H\cdot\mathbf{r}_{t})
-(\mathbf{r}_s\cdot\mathbf{r}_t)(\mathbf{r}_{t}*H\cdot\mathbf{r}_{s})
\cr &&
+(\mathbf{r}_t\cdot\mathbf{r}_t)(\mathbf{r}_{s}*H\cdot\mathbf{r}_{s})
-(\mathbf{r}_t\cdot\mathbf{r}_s)(\mathbf{r}_{s}*H\cdot\mathbf{r}_{t}))
\cr &=&\frac{-1}{|\nabla F|}((\mathbf{r}_s\times
(\mathbf{r}_t*H))\cdot(\mathbf{r}_s\times
\mathbf{r}_t)+(\mathbf{r}_t\times
(\mathbf{r}_s*H))\cdot(\mathbf{r}_t\times \mathbf{r}_s)) \cr
&=&\frac{1}{|\nabla F|}(\mathbf{r}_t\times
(\mathbf{r}_s*H)-\mathbf{r}_s\times
(\mathbf{r}_t*H))\cdot(\mathbf{r}_s\times \mathbf{r}_t) \cr
&=&\frac{1}{|\nabla F|}((\mathbf{r}_s\times
\mathbf{r}_t)*H-Trace(H)(\mathbf{r}_s\times
\mathbf{r}_t))\cdot(\mathbf{r}_s\times \mathbf{r}_t) \cr
&=&\frac{1}{|\nabla F|}((\mathbf{r}_s\times
\mathbf{r}_t)H(\mathbf{r}_s\times
\mathbf{r}_t)^T-Trace(H)|\mathbf{r}_s\times \mathbf{r}_t|^2).
\end{eqnarray*}
Hence, the theorem can be proved by substituting the above results
into formulas (\ref{eq:qlg}) and (\ref{eq:qlm}).
\end{pf*}

\begin{example}{\rm
$F(x,y,z)=(x-y)(x^2+y^2+z^2-2Rx)$. The surface $F=0$ is composed of
a plane $x=y$ and a sphere with a radius of $R$ crossing at the
origin point $P=(0,0,0)$, as shown in Fig.~\ref{fig:surface1}. It is
clear that $P$ is singular point of the surface since $\nabla
P(0,0,0)=0$. There are two tangent planes with normal vectors
$(1,-1,0)$ and $(1,0,0)$ corresponding to the plane and the sphere
at $P$, respectively.

By Algorithm \ref{Alg:Sql} and the quadratic parametric surface
Eq.~(\ref{eq:PS}), we have
$$F(x,y,z)=C_{1,0}s+C_{0,1} t+C_{2,0} s^2+C_{1,1} st +C_{0,2}
t^2+\cdots,$$ where $C_{1,0}=C_{0,1}=0$, $C_{2,0}=2a_1(b_1-a_1)$,
$C_{1,1}=2a_1(b_2-a_2)R+2a_2(b_1-a_1)R$, $C_{0,2}=2a_2(b_2-a_2)R$.
The origin point is a 2-fold point as $C_1=0,\ C_2\neq 0$.

Let $C_{2,0}=C_{1,1}=C_{0,2}=0$, we obtain two solutions $a_1=b_1,\
a_2=b_2$ and $a_1=0,\ a_2=0$.
\begin{enumerate}[\mbox{Case }1:]
\item If $a_1=b_1,\ a_2=b_2$, then the normal vector of the tangent plane
is
$\mathbf{r}_s\times\mathbf{r}_t=(a_1c_2-a_2c_1,-a_1c_2+a_2c_1,0)$,
equivalent to $(1,-1,0)$. The coefficients of $s^3,\ s^2t,\ st^2,\
t^3$ are
$$C_{3,0}=a_1(b_3-a_3)R,\ C_{2,1}=2a_1(b_4-a_4)R+a_2(b_3-a_3)R,$$
$$C_{1,2}=2a_2(b_4-a_4)R+a_1(b_5-a_5)R,\ C_{0,3}=a_2(b_5-a_5)R.$$
Let $C_{3,0}=C_{2,1}=C_{1,2}=C_{0,3}=0$, then by Eq.~(\ref{eq:qlg})
and Eq.~(\ref{eq:qlm}), we obtain the Gaussian and mean curvatures
at $P$ are $K_G=0,\ K_M=0$. This case corresponds to the plane
branch of the surface at the origin.

\item If $a_1=0,\ a_2=0$, then the normal vector of the tangent plane
is $\mathbf{r}_s\times\mathbf{r}_t=(b_1c_2-b_2c_1,0,0)$, equivalent
to $(1,0,0)$. The coefficients of $s^3,\ s^2t,\ st^2,\ t^3$ are
$$C_{3,0}=b_1(a_3R-b_1^2-c_1^2),\ C_{2,1}=b_2(a_3R-b_1^2-c_1^2)+2b_1(a_4R-b_1b_2-c_1c_2),$$
$$C_{1,2}=b_1(a_5R-b_2^2-c_2^2)+2b_2(a_4R-b_1b_2-c_1c_2),\ C_{0,3}=b_2(a_5R-b_2^2-c_2^2).$$
Let $C_{3,0}=C_{2,1}=C_{1,2}=C_{0,3}=0$, then by Eq.~(\ref{eq:qlg})
and Eq.~(\ref{eq:qlm}), we obtain the Gaussian and mean curvatures
at $P$ are $K_G=1/R^2,\ K_M=1/R$. This case corresponds to the
sphere branch of the surface at the origin.
\end{enumerate}
\begin{figure}[!h]
\centering
\includegraphics[width=0.35\textwidth]{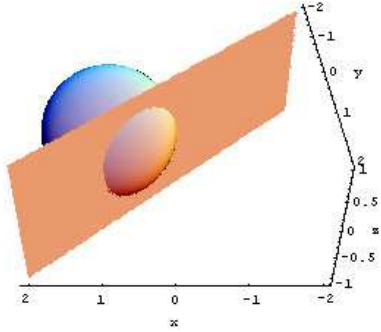}
\caption{The surface
$(x-y)(x^2+y^2+z^2-2Rx)=0$.}\label{fig:surface1}
\end{figure}
}\end{example}

\begin{example}\label{Ex:S2}{\rm
$F(x,y,z)=x^4+y^2+yz^2-z^2$. The surface $F=0$ is composed of two
components intersecting at the origin point $P=(0,0,0)$, as shown in
Fig.~\ref{fig:surface2}. It is clear that $P$ is singular point of
the surface since $\nabla P(0,0,0)=0$. By Algorithm \ref{Alg:Sql}
and the quadratic parametric surface Eq.~(\ref{eq:PS}), we have
$$F(x,y,z)=C_{1,0}s+C_{0,1} t+C_{2,0} s^2+C_{1,1} st +C_{0,2}
t^2+\cdots,$$ where $C_{1,0}=C_{0,1}=0$, $C_{2,0}=b_1^2-c_1^2$,
$C_{1,1}=2(b_1b_2-c_1c_2)$, $C_{0,2}=b_2^2-c_2^2$. Hence, the origin
point is a 2-fold point.

Let $C_{2,0}=C_{1,1}=C_{0,2}=0$, we obtain two solutions $b_1=c_1,\
b_2=c_2$ and $b_1=-c_1,\ b_2=-c_2$.
\begin{enumerate}[\mbox{Case }1:]
\item If $b_1=c_1,\
b_2=c_2$, then the normal vector of the tangent plane is
$\mathbf{r}_s\times\mathbf{r}_t=(0,-a_1b_2+a_2c_1,a_1b_2-a_2c_1)$,
equivalent to $(0,1,-1)$. The coefficients of $s^3,\ s^2t,\ st^2,\
t^3$ are
$$C_{3,0}=c_1(c_1^2+b_3-c_3),\ C_{2,1}=b_2(b_3-c_3)+3b_2c_1^2+2c_1(b_4-c_4),$$
$$C_{1,2}=2b_2(b_4-c_4)+3b_2^2c_1+c_1(b_5-c_5),\ C_{0,3}=b_2(b_2^2+b_5-c_5).$$
Let $C_{3,0}=C_{2,1}=C_{1,2}=C_{0,3}=0$, then by Eq.~(\ref{eq:qlg})
and Eq.~(\ref{eq:qlm}), we obtain the Gaussian and mean curvatures
at $P$ are $K_G=0,\ |K_M|=\sqrt{2}/8$.

\item If $b_1=-c_1,\
b_2=-c_2$, then the normal vector of the tangent plane is
$\mathbf{r}_s\times\mathbf{r}_t=(0,a_1b_2+a_2c_1,a_1b_2+a_2c_1)$,
equivalent to $(0,1,1)$. The coefficients of $s^3,\ s^2t,\ st^2,\
t^3$ are
$$C_{3,0}=-c_1(c_1^2+b_3+c_3),\ C_{2,1}=b_2(b_3+c_3)+3b_2c_1^2-2c_1(b_4+c_4),$$
$$C_{1,2}=2b_2(b_4+c_4)-3b_2^2c_1-c_1(b_5+c_5),\ C_{0,3}=b_2(b_2^2+b_5+c_5).$$
Let $C_{3,0}=C_{2,1}=C_{1,2}=C_{0,3}=0$, then by Eq.~(\ref{eq:qlg})
and Eq.~(\ref{eq:qlm}), we obtain the Gaussian and mean curvatures
at $P$ are $K_G=0,\ |K_M|=\sqrt{2}/8$.
\end{enumerate}
The above two cases are corresponding with the two components of the
surface $F(x,y,z)=0$ at the singular point $P$, as shown in
Fig.~\ref{fig:surface2}.

\begin{figure}[!h]
\centering
\includegraphics[width=0.35\textwidth]{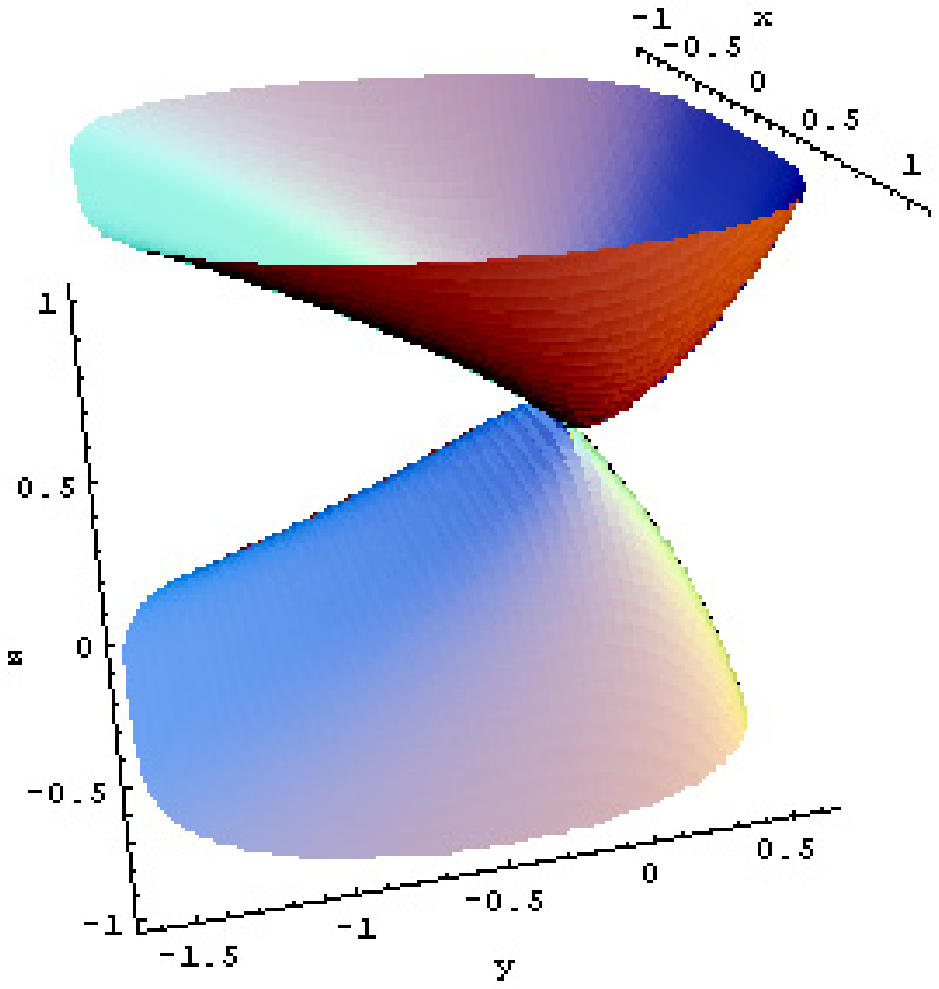}\hspace{5ex}
\includegraphics[width=0.35\textwidth]{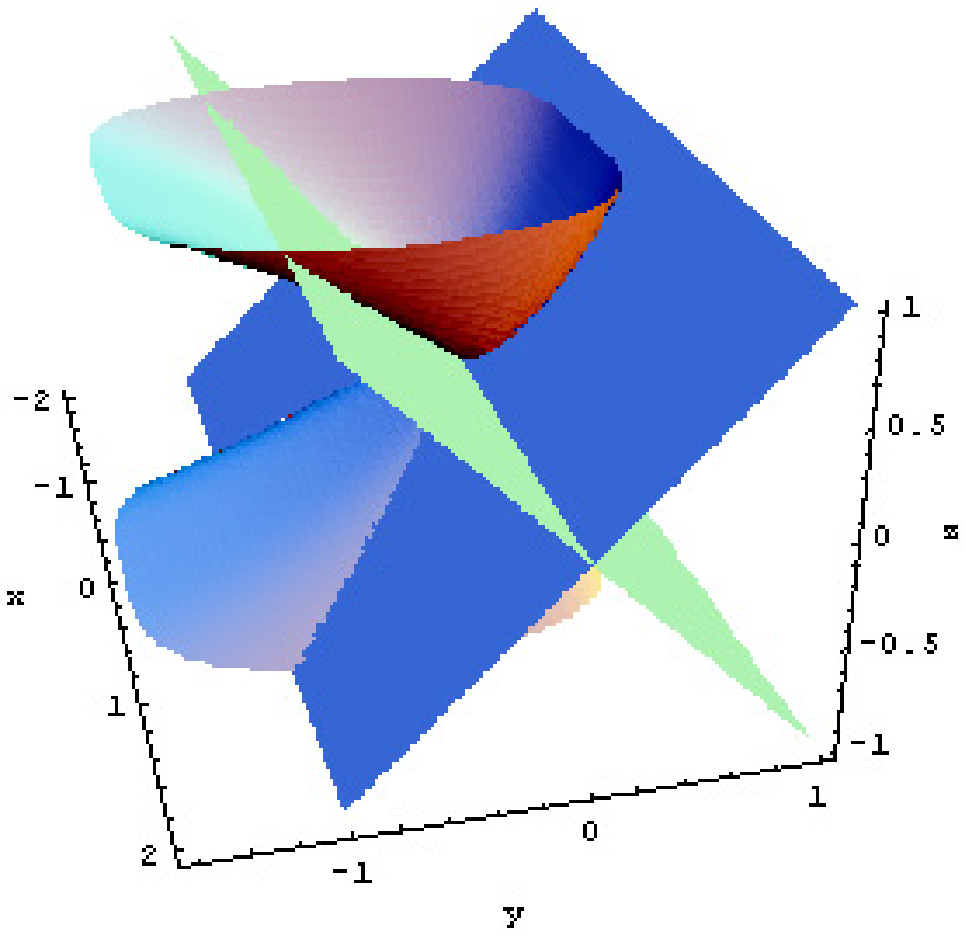}
\caption{The surface $x^4+y^2+yz^2-z^2=0$ and its two tangent plane
at the singular point $P=(0,0,0)$.}\label{fig:surface2}
\end{figure}
}\end{example}

\section{Curvatures and torsions at the singular points of space algebraic curves}
For given a space algebraic curve $\mathcal{C}$ defined by the
intersection of two algebraic surfaces $\mathcal{S}_1: F(x,y,z)=0$
and $\mathcal{S}_2: G(x,y,z)=0$, i.e.,
\begin{equation}\label{eq:spacecurve}
\mathcal{C}: \left\{\begin{array}{ll}F(x,y,z)=0,\\
G(x,y,z)=0.\end{array}\right. \end{equation} A point $P$ is called a
singular point of $\mathcal{C}$, if where $\nabla F\times \nabla G
=0$. It includes three cases: $\nabla F=0$, $\nabla G=0$, or $\nabla
F=\lambda\nabla G$ for a constant $\lambda\neq 0$.

Firstly, we discuss the multiplicity of the singularity. Similar to
the case of plane curves, we can define the multiplicity of the
singular point by the number of tangents (counting multiplicity) of
the curve at the singular point. Note that each tangent line to the
space curve $\mathcal{C}$ is determined by the intersection of a
pair of tangent planes to the two surfaces $\mathcal{S}_1$ and
$\mathcal{S}_2$ respectively. Suppose that the multiplicities of $P$
are $r_1$ and $r_2$ corresponding to surfaces $\mathcal{S}_1$ and
$\mathcal{S}_2$ respectively. It means that $\mathcal{S}_1$ has
$r_1$ tangent planes and $\mathcal{S}_2$ has $r_2$ tangent planes at
$P$. If the $r_1$ tangent planes are different to the $r_2$ tangent
planes, then the curve $\mathcal{C}$ has $r_1r_2$ tangent lines at
$P$. However, we remark that the multiplicity of the singularity of
the space curve is more complex than the case of the plane curve,
since the space curve is determined by two equations. We present the
computation of the multiplicity of the singularity at a point as
follows.

Let $\Gamma$ be a cubic parametric curve through the point $P$ of
$\mathcal{C}$, defined by the following equations
\begin{equation} \label{eq:PC3}
\mathbf{r}(t)=\left\{
\begin{aligned}
x(t)&=a+a_1 t+\frac{1}{2!}a_2 t^2+\frac{1}{3!}a_3 t^3, \\
y(t)&=b+b_1 t+\frac{1}{2!}b_2 t^2+\frac{1}{3!}b_3 t^3, \\
z(t)&=c+c_1 t+\frac{1}{2!}c_2 t^2+\frac{1}{3!}c_3 t^3.
\end{aligned} \right.
\end{equation}
Then $\mathbf{r}(0)=(a,b,c)$, $\mathbf{r}'(0)=(a_1,b_1,c_1)$,
$\mathbf{r}''(0)=(a_2,b_2,c_2)$, $\mathbf{r}'''(0)=(a_3,b_3,c_3)$.
The curvature and torsion at $P=(a,b,c)$ of $\Gamma$ are
\begin{eqnarray}
k&=&\frac{|\mathbf{r}'\times \mathbf{r}''|}{|\mathbf{r}'|^3}, \label{eq:ql2}\\
\tau&=&\frac{\det(\mathbf{r}'\ \mathbf{r}''\
\mathbf{r}''')}{|\mathbf{r}'\times \mathbf{r}''|^2}.\label{eq:nl}
\end{eqnarray}

The intersections of $\Gamma$ and $\mathcal{C}$ are corresponding to
the roots of the equations $F(x(t),y(t),z(t))=0$ and
$G(x(t),y(t),z(t))=0$. Expand them by the power of $t$ and notice
that $F(a,b,c)=G(a,b,c)=0$, then we have
\begin{eqnarray}\label{eq:FGt}
\left\{\begin{array}{ll}F(x,y,z)=C_{1}t+C_{2}t^2+C_{3}t^3+\cdots,\\
G(x,y,z)=D_{1}t+D_{2}t^2+D_{3}t^3+\cdots.\end{array}\right.
\end{eqnarray}
In order that the above two equations have common root of $t$, it is
necessary that they have sufficient same coefficients $C_i=D_i,\
i=1,2,\ldots$.

Suppose that $l$ is the minimal natural number such that $C_1=D_1=0,
\ldots, C_{l-1}=D_{l-1}=0$, $C_l=D_l\neq 0$. If there are $r$
solutions of $(a_1,b_1,c_1)$ (counting multiplicity) by $C_l=D_l=0$
(if needed, let $C_{l+1}=D_{l+1},\ldots$), then we obtain the $r$
tangents to the curve at $P$ and the multiplicity of the singularity
of $P$ is $r$.

Further, the curvature and torsion at $P$ along every tangent to
$\mathcal{C}$ is determined by the cubic curve $\Gamma$, which has
intersection with $\mathcal{C}$ at $P$ with multiplicity of much
more than $r$. We conclude the procedure as the following algorithm.
\begin{algor}[Curvatures at the singular point of a space algebraic
curve]\label{Alg:ql3}

For given a space algebraic curve $\mathcal{C}$ defined by
Eq.~(\ref{eq:spacecurve}), $P=(a,b,c)$ is an $r$-fold point of
$\mathcal{C}$. Let $\Gamma$ be a parametric cubic curve through the
point $P$ defined by Eq.~(\ref{eq:PC3}).
\begin{enumerate}[1)]
\item By Eq.~(\ref{eq:FGt}), if $l(\geq 1)$ is the minimal natural number such that
$C_1=\cdots=C_{l-1}=D_1=\cdots=D_{l-1}=0,C_l=D_l\neq0$, then the $r$
tangents $(a_1,b_1,c_1)$ of $\mathcal{C}$ at $P$ are determined by
equations $C_l=D_l=0$ (if needed, let $C_{l+1}=D_{l+1},\ldots$).

\item Substitute the above $a_1,\ b_1,\ c_1$ into Eq.~(\ref{eq:FGt}). If
$p(>l)$ is the minimal natural number such that
$C_1=\cdots=C_{p-1}=D_1=\cdots=D_{p-1}=0,C_p=D_p\neq 0$, then the
curvature of $\mathcal{C}$ at $P$ along the tangent $(a_1,b_1,c_1)$
equals to the curvature of $\Gamma$ at $P$, which determined by
equations $C_p=D_p=0$ (if needed, let $C_{p+1}=D_{p+1},\ldots$). The
curvature is computed by Eq.~(\ref{eq:ql2}).

\item Substitute the above $a_2,\ b_2,\ c_2$ into Eq.~(\ref{eq:FGt}). If
$q(>p)$ is the minimal natural number such that
$C_1=\cdots=C_{q-1}=D_1=\cdots=D_{q-1}=0,C_q=D_q\neq0$, then the
torsion of $\mathcal{C}$ at $P$ along the tangent $(a_1,b_1,c_1)$
equals to the torsion of $\Gamma$ at $P$, which determined by
equations $C_q=D_q=0$ (if needed, let $C_{q+1}=D_{q+1},\ldots$), The
torsion is computed by Eq.~(\ref{eq:nl}).
\end{enumerate}
\end{algor}

In order to verify the correctness of Algorithm \ref{Alg:ql3}, we
prove that the curvature and torsion obtained by
Algorithm~\ref{Alg:ql3} is same to that by formula in
\cite{Goldman2005} when the point $P$ is regular.

\begin{theorem}
If $P$ is a regular point of a space algebraic curve $\mathcal{C}$,
the curvature and torsion at $P$ obtained by Algorithm \ref{Alg:ql3}
are
\begin{eqnarray*}
k&=&\frac{|\mathbf{r}'\times \mathbf{r}''|}{|\mathbf{r}'|^3}
=\frac{|((\nabla F\times \nabla G)\cdot\nabla(\nabla F\times \nabla G))\times(\nabla F\times \nabla G)|}{|\nabla F\times \nabla G|^3}, \\
\tau&=&\frac{\det(\mathbf{r}'\ \mathbf{r}''\
\mathbf{r}''')}{|\mathbf{r}'\times \mathbf{r}''|^2}=\frac{\det(T^*\
T^{**}\ T^{***})}{|T^*\times T^{**}|^2}.
\end{eqnarray*}
where $T^*,T^{**},T^{***}$ are defined as \cite{Goldman2005},
$$T^*=\nabla F\times \nabla G,\ T^{**}=(\nabla F\times \nabla
G)\cdot\nabla(\nabla F\times \nabla G),$$
$$T^{***}=(\nabla F\times \nabla G)\cdot\nabla(\nabla(\nabla F\times \nabla G))\cdot(\nabla F\times \nabla G)
+(\nabla F\times \nabla G)\cdot\nabla(\nabla F\times \nabla
G)\cdot\nabla(\nabla F\times \nabla G).$$
\end{theorem}
\begin{pf*}{Proof}
The multiplicity of $P$ is $r=1$, $\nabla F\neq 0$, $\nabla G\neq 0$
and $\nabla F\times \nabla G\neq 0$, then
\begin{eqnarray}\label{eq:C2}
\begin{array}{ll}
C_{1}&=F_xa_1+F_yb_1+F_zc_1=\mathbf{r}'\cdot\nabla F,\\
C_{2}&=\frac{1}{2}(\mathbf{r}''\cdot\nabla F+\mathbf{r}'\cdot H(f)\cdot\mathbf{r}'),\\
C_{3}&=\frac{1}{6}(\mathbf{r}'''\cdot\nabla F+3\mathbf{r}''\cdot H(f)\cdot\mathbf{r}'+(\mathbf{r}'\cdot \nabla H(f)\cdot\mathbf{r}')\mathbf{r}'),\\
D_{1}&=G_xa_1+G_yb_1+G_zc_1=\mathbf{r}'\cdot\nabla G,\\
D_{2}&=\frac{1}{2}(\mathbf{r}''\cdot\nabla G+\mathbf{r}'\cdot H(G)\cdot\mathbf{r}'),\\
D_{3}&=\frac{1}{6}(\mathbf{r}'''\cdot\nabla G+3\mathbf{r}''\cdot H(G)\cdot\mathbf{r}'+(\mathbf{r}'\cdot \nabla H(G)\cdot\mathbf{r}')\mathbf{r}'),\\
\end{array}
\end{eqnarray}
where $\mathbf{r}'\cdot \nabla F$ means the inner product of two
vectors $\mathbf{r}'$ and $\nabla F$, $\mathbf{r}'\cdot
H(f)\cdot\mathbf{r}'$ means product of row vector $\mathbf{r}'$,
matrix $H(f)$ and column vector $\mathbf{r}'^T$, and $\nabla$
applied to a matrix such as $H(F)$ means apply $\nabla$ to each
column vector of the matrix to generate a list of three consecutive
matrices \cite{Goldman2005}.

By $C_{1}=D_{1}=0$, we obtain $\mathbf{r}'=\lambda(\nabla F\times
\nabla G)$, for a constant $\lambda\neq
 0$. Let $\lambda=1$ for convenience, then $\mathbf{r}'=\nabla F\times \nabla
 G$.

By $C_{2}=D_{2}=0$, we obtain $\mathbf{r}''\cdot\nabla
F=-\mathbf{r}'\cdot H(F)\cdot \mathbf{r}'$ and
$\mathbf{r}''\cdot\nabla G=-\mathbf{r}'\cdot H(G)\cdot \mathbf{r}'$.
Then, $$\mathbf{r}'\times \mathbf{r}''=(\nabla F\times \nabla
G)\times \mathbf{r}''=(\nabla F\cdot \mathbf{r}'')\nabla G-(\nabla
G\cdot \mathbf{r}'')\nabla F=(\mathbf{r}'\cdot H(G)\cdot
\mathbf{r}')\nabla F-(\mathbf{r}'\cdot H(F)\cdot \mathbf{r}')\nabla
G.$$ By computation, it is easy to verify that $$(\mathbf{r}'\cdot
H(G)\cdot \mathbf{r}')\nabla F-(\mathbf{r}'\cdot H(F)\cdot
\mathbf{r}')\nabla G=(\nabla F\times \nabla G)\times((\nabla F\times
\nabla G)\cdot\nabla(\nabla F\times \nabla G)).$$ That is
$$(\mathbf{r}'\cdot\nabla(\nabla F\times \nabla G))\times \mathbf{r}'= \mathbf{r}''\times \mathbf{r}'.$$

Note that
\begin{eqnarray*}&&\nabla(\nabla F\times \nabla G)\cr
&=&\nabla\{\det\binom{F_y\ F_z}{G_y\ G_z},\ \det\binom{F_z\
F_x}{G_z\ G_x},\ \det\binom{F_x\ F_y}{G_x\ G_y}\}\cr &=& H(F)\times
\nabla G-H(G)\times \nabla F\cr &=&\{\nabla F_x\times\nabla G,\
\nabla F_y\times\nabla G,\ \nabla F_z\times\nabla G\}-\{\nabla
G_x\times\nabla F,\ \nabla G_y\times\nabla F,\ \nabla
G_z\times\nabla F\},
\end{eqnarray*}
then \begin{eqnarray*}&&\nabla(\nabla F\times \nabla G)\cdot\nabla F
\cr &=& \{(\nabla F_x\times\nabla G)\cdot \nabla F,\ (\nabla
F_y\times\nabla G)\cdot \nabla F,\ (\nabla F_z\times\nabla G)\cdot
\nabla F\} \cr &=& -(\nabla F\times \nabla G)\{\nabla F_x,\ \nabla
F_y,\ \nabla F_z\}=-(\nabla F\times \nabla G)\cdot H(F) \cr
&=&-\mathbf{r}'\cdot H(F)=-H(F)\cdot\mathbf{r}',
\end{eqnarray*}we have
\begin{eqnarray*}&&
(\mathbf{r}'\cdot\nabla(\nabla F\times \nabla G))\cdot \nabla F
=\mathbf{r}'\cdot(\nabla(\nabla F\times \nabla G)\cdot \nabla F)
=-\mathbf{r}'\cdot H(F)\cdot\mathbf{r}'=\mathbf{r}''\cdot \nabla F,
\end{eqnarray*}
similarly, $(\mathbf{r}'\cdot\nabla(\nabla F\times \nabla G))\cdot
\nabla G=\mathbf{r}''\cdot \nabla G$. Hence, we obtain
$\mathbf{r}''=\mathbf{r}'\cdot\nabla(\nabla F\times \nabla G)$.

By $C_3=D_3=0$, we have
\begin{eqnarray*}
&&(\mathbf{r}',\mathbf{r}'',\mathbf{r}''')=(\mathbf{r}'''\times
\mathbf{r}')\cdot \mathbf{r}''=-(\mathbf{r}'\times
\mathbf{r}''')\cdot \mathbf{r}''=-((\nabla F\times \nabla G)\times
\mathbf{r}''')\cdot \mathbf{r}'' \cr &=& -((\nabla F\cdot
\mathbf{r}''')\nabla G-(\nabla G\cdot \mathbf{r}''')\nabla F)\cdot
\mathbf{r}'' \cr &=& (\mathbf{r}'''\cdot \nabla G)(\mathbf{r}''\cdot
\nabla F)-(\mathbf{r}'''\cdot \nabla F)(\mathbf{r}''\cdot \nabla G)
\cr &=& (3\mathbf{r}''\cdot H(G)\cdot \mathbf{r}'+(\mathbf{r}'\cdot
\nabla H(G)\cdot \mathbf{r}')\cdot \mathbf{r}')(\mathbf{r}'\cdot
H(F)\cdot \mathbf{r}') \cr && -(3\mathbf{r}''\cdot H(F)\cdot
\mathbf{r}'+(\mathbf{r}'\cdot \nabla H(F)\cdot \mathbf{r}')\cdot
\mathbf{r}')(\mathbf{r}'\cdot H(G)\cdot \mathbf{r}').
\end{eqnarray*}
By computation, it equals to $(T^*,T^{**},T^{***})$, where
$T^*=\mathbf{r}'$, $T^{**}=\mathbf{r}''$ and
$$T^{***}=(\nabla F\times \nabla G)\cdot\nabla(\nabla(\nabla F\times \nabla G))\cdot(\nabla F\times \nabla G)
+(\nabla F\times \nabla G)\cdot\nabla(\nabla F\times \nabla
G)\cdot\nabla(\nabla F\times \nabla G)$$ according to
\cite{Goldman2005}. So we prove that the curvature
Eq.~(\ref{eq:ql2}) and torsion Eq.~(\ref{eq:nl}) at $P$ of the
parametric curve $\Gamma$ are equivalence to that presented in
\cite{Goldman2005}.
\end{pf*}

\begin{figure}[!h]
\centering
\begin{minipage}[t]{0.45\textwidth}
\includegraphics[width=0.8\textwidth]{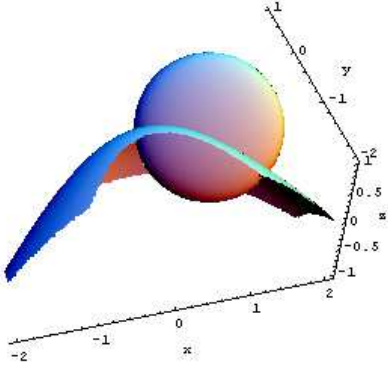}
\caption{The figure of the curve $\{x^2+y^2+z^2-2Rx=0\}\cap
\{x^2+2y-yz+z=0\}$.}\label{fig:S-curve-1}
\end{minipage}\hspace{3ex}
\begin{minipage}[t]{0.45\textwidth}
\includegraphics[width=0.8\textwidth]{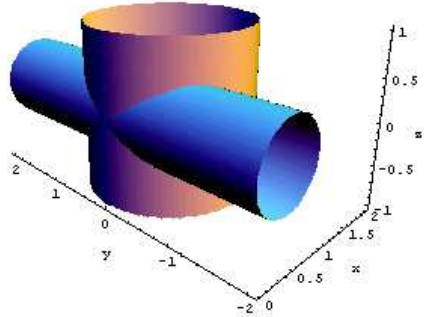}
\caption{The figure of the curve
$\{x^2+y^2-2R_1x=0\}\cap\{x^2+z^2-2R_2x=0\}$.}\label{fig:S-curve-2}
\end{minipage}
\end{figure}

Lastly, we apply Algorithm \ref{Alg:ql3} to the following examples.
\begin{example}
$F(x,y,z)=x^2+y^2+z^2-2Rx$, $G(x,y,z)=x^2+2y-yz+z$. $\nabla
F(0,0,0)=(-2R,0,0)$, $\nabla G(0,0,0)=(0,2,1)$, $\nabla F\times
\nabla G\neq 0$, then the origin $P=(0,0,0)$ is a regular point (as
shown in Fig.~\ref{fig:S-curve-1}).

By Algorithm \ref{Alg:ql3} and a parametric cubic curve
(\ref{eq:PC3}), we have
\begin{eqnarray*}
\left\{\begin{array}{ll}F(x,y,z)=C_{1}t+C_{2}t^2+C_{3}t^3+\cdots,\\
G(x,y,z)=D_{1}t+D_{2}t^2+D_{3}t^3+\cdots,\end{array}\right.
\end{eqnarray*}
where, $C_1=-a_1 R$, $D_1=2b_1+c_1$.
\begin{enumerate}[1)]
\item By $C_1=D_1=0$, we obtain $a_1=0,\ c_1=-2b_1$,
then the tangent vector to the curve at $P$ is $(0,1,-2)$, and
$C_2=5b_1^2-a_2R$, $D_2=2b_1^2+b_2+c_2/2$.

\item By $C_2=D_2=0$, we obtain $a_2=5b_1^2/R,\ c_2=-2(2b_1^2+b_2)$,
By Eq.~(\ref{eq:ql2}), the curvature of the curve at $P$ is
$\sqrt{\frac{125+16R^2}{125R^2}}$, and $C_3=8b_1^3+5b_1b_2-a_3R/3,\
D_3=2b_1^3+2b_1b_2+b_3/3+c_3/6$.

\item By $C_3=D_3=0$, we obtain $a_3=3(8b_1^3+5b_1b_2)/R,\ c_3=-2(6b_1^3+6b_1b_2+b_3)$,
By Eq.~(\ref{eq:nl}), the torsion of the curve at $P$ is
$-\frac{36R}{125+16R^2}$.
\end{enumerate}
\end{example}

\begin{example}
$F(x,y,z)=x^2+y^2-2R_1x$, $G(x,y,z)=x^2+z^2-2R_2x$. $\nabla
F(0,0,0)=(-2R_1,0,0)$, $\nabla G(0,0,0)=(-2R_2,0,0)$, $\nabla
F\times \nabla G= 0$, then the origin $P=(0,0,0)$ is a singular
point of the curve (as shown in Fig.~\ref{fig:S-curve-2}).

By Algorithm \ref{Alg:ql3} and a parametric cubic curve
(\ref{eq:PC3}), we have
\begin{eqnarray*}
\left\{\begin{array}{ll}F(x,y,z)=C_{1}t+C_{2}t^2+C_{3}t^3+\cdots,\\
G(x,y,z)=D_{1}t+D_{2}t^2+D_{3}t^3+\cdots,\end{array}\right.
\end{eqnarray*}
where, $C_1=-2a_1 R_1,\ D_1=-2a_1 R_2,\ C_2=a_1^2+b_1^2-a_2R_1,\
D_2=a_1^2+c_1^2-a_2R_2$.

Case 1: If $R_1=R_2=R$, then $C_1=D_1$. By $C_2=D_2$, we obtain
$b_1=\pm c_1$. Then let $C_1=D_1=0$, we have $a_1=0$. Hence, we
obtain two tangents $(0,1,1)$ and $(0,1,-1)$, and the point $P$ is a
2-fold point of the curve.
\begin{enumerate}[1)]
\item If $a_1=0,\ b_1=c_1$,
the tangent to the curve at $P$ is $(0,1,1)$, and
$C_2=D_2=b_1^2-a_2R$, $C_3=b_1b_2-a_3R/3$, $D_3=b_1c_2-a_3R/3$.

By $C_2=D_2=0,\ C_3=D_3$, we obtain $a_2=b_1^2/R,\ b_2=c_2$. By
Eq.~(\ref{eq:ql2}), the curvature of the curve at $P$ is
$\frac{1}{2R}$. Then let $C_3=D_3=0,\ C_4=D_4$, by
Eq.~(\ref{eq:nl}), the torsion of the curve at $P$ is $0$.

\item Similarly, we obtain the curvature and torsion of the curve at $P$
along the tangent $(0,1,-1)$ are $\frac{1}{2R}$ and $0$.
\end{enumerate}

Case 2: If $R_1\neq R_2$. By $C_1=D_1,\ C_2=D_2$, we obtain $a_1=0,\
a_2=(b_1^2-c_1^2)/(R_1-R_2)$, and $C_1=D_1=0,\
C_2=D_2=(b_1^2R_2-c_1^2R_1)/(R_2-R_1)$. By $C_2=D_2=0$, then we
obtain two solutions $a_1=0,\ c_1=\pm b_1\sqrt{R_2/R_1}$, and the
point $P$ is a 2-fold point of the curve.
\begin{enumerate}[1)]
\item If $a_1=0,\ c_1=b_1\sqrt{R_2/R_1}$,
the tangent to the curve at $P$ is $(0,1,\sqrt{R_2/R_1})$, and
$C_3=b_1b_2-a_3R_1/3$, $D_3=b_1c_2\sqrt{R_2/R_1}-a_3R_2/3$.

By $C_3=D_3=0$, we obtain $b_2=\frac{a_3R_1}{3b_1},\
c_2=\frac{a_3\sqrt{R_1R_2}}{3b_1}$, By Eq.~(\ref{eq:ql2}), the
curvature of the curve at $P$ is $1/(R_1+R_2)$. Then let
$C_4=D_4=0$, by Eq.~(\ref{eq:nl}), the torsion of the curve at $P$
is $\frac{3(R_1-R_2)}{4\sqrt{R_1R_2}(R_1+R_2)}$.

\item Similarly, we obtain the curvature and torsion of the curve at $P$
along the tangent $(0,1,-\sqrt{R_2/R_1})$ are $1/(R_1+R_2)$ and
$\frac{3(R_2-R_1)}{4\sqrt{R_1R_2}(R_1+R_2)}$.
\end{enumerate}
\end{example}

\begin{figure}[!h]
\centering
\begin{minipage}[t]{0.45\textwidth}
\includegraphics[width=0.8\textwidth]{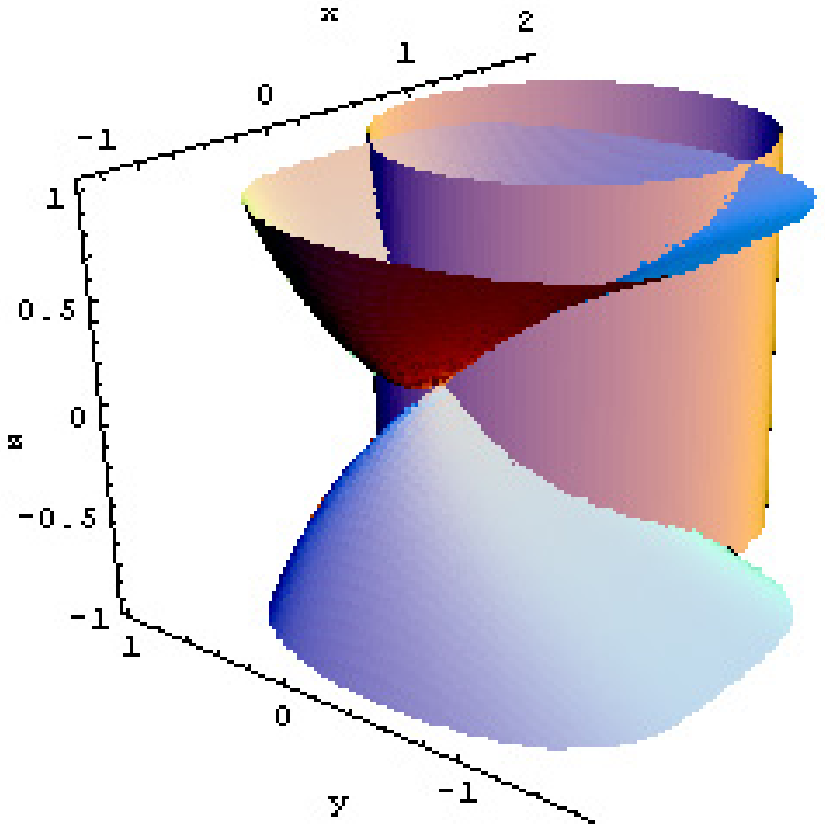}
\caption{The figure of the curve $\{x^4+y^2+yz^2-z^2=0\}\cap
\{x^2+y^2-2Rx=0\}$.}\label{fig:S-curve-3}
\end{minipage}\hspace{3ex}
\begin{minipage}[t]{0.45\textwidth}
\includegraphics[width=0.8\textwidth]{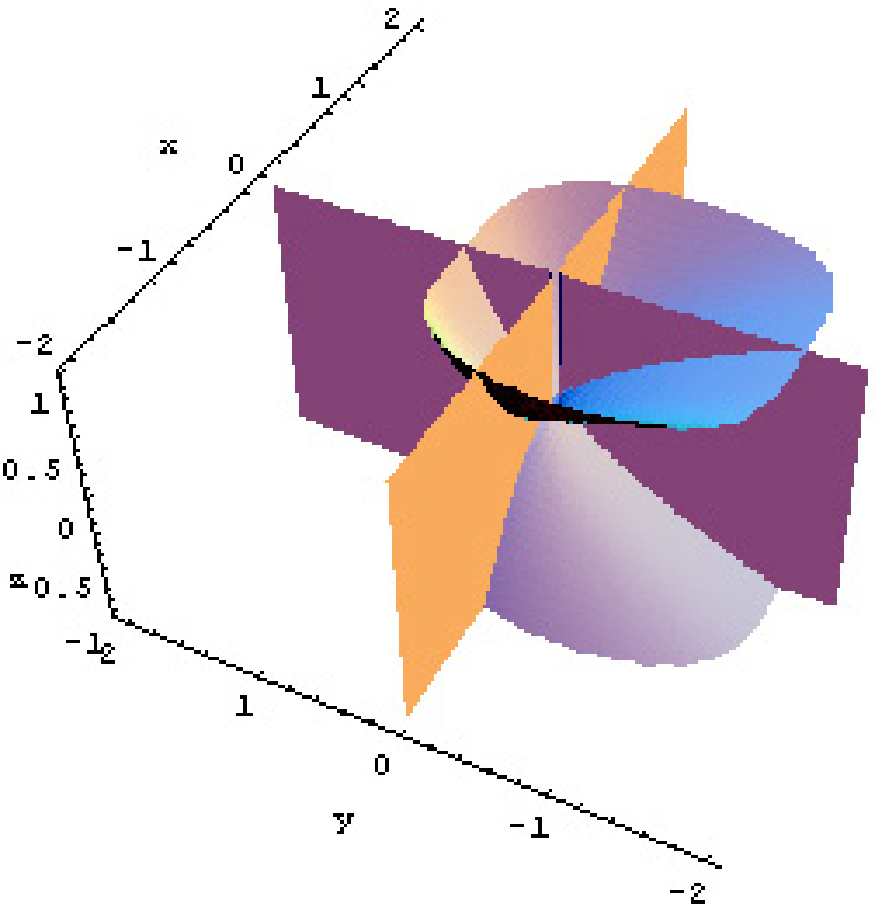}
\caption{The figure of the curve $\{x^4+y^2+yz^2-z^2=0\}\cap
\{xy=0\}$.}\label{fig:S-curve-4}
\end{minipage}
\end{figure}

\begin{example}
$F(x,y,z)=x^4+y^2+yz^2-z^2$, $G(x,y,z)=x^2+y^2-2Rx$. $\nabla
F(0,0,0)=(0,0,0)$, $\nabla G(0,0,0)=(-2R,0,0)$. As discussed in
Example \ref{Ex:S2}, the origin $P=(0,0,0)$ is a 2-fold point of the
surface $F(x,y,z)=0$, and is a simple point of the surface
$G(x,y,z)=0$. Then $P$ is a 2-fold point of the curve since the
tangent planes of the two surfaces are different (as shown in
Fig.~\ref{fig:S-curve-3}).

By Algorithm \ref{Alg:ql3} and a parametric cubic curve
(\ref{eq:PC3}), we have
\begin{eqnarray*}
\left\{\begin{array}{ll}F(x,y,z)=C_{1}t+C_{2}t^2+C_{3}t^3+\cdots,\\
G(x,y,z)=D_{1}t+D_{2}t^2+D_{3}t^3+\cdots,\end{array}\right.
\end{eqnarray*}
where, $C_1=0,\ D_1=-2a_1 R,\ C_2=b_1^2-c_1^2,\
D_2=a_1^2+b_1^2-a_2R$.

By $C_1=D_1,\ C_2=D_2$, we obtain $a_1=0,\ a_2=c_1^2/R$, and
$C_1=D_1=0,\ C_2=D_2=b_1^2-c_1^2$. By $C_2=D_2=0$, then $a_1=0,\
c_1=\pm b_1$.
\begin{enumerate}[1)]
\item If $a_1=0,\ c_1=b_1$,
then the tangent to the curve at $P$ is $(0,1,1)$, and
$C_3=b_1^3+b_1b_2-b_1c_2$, $D_3=b_1c_2-a_3R/3$.

Let $C_3=D_3=0$, we obtain $a_3=3b_1b_2/R,\ c_2=b_1^2+b_2$. By
Eq.~(\ref{eq:ql2}), the curvature of the curve at $P$ is
$\frac{\sqrt{2+R^2}}{2\sqrt{2}R}$. Then let $C_4=D_4=0$, by
Eq.~(\ref{eq:nl}), the torsion of the curve at $P$ is
$-\frac{9R}{4(2+R^2)}$.

\item Similarly, we obtain the curvature and torsion of the curve at $P$
along the tangent $(0,1,-1)$ are $\frac{\sqrt{2+R^2}}{2\sqrt{2}R}$
and $\frac{9R}{4(2+R^2)}$.
\end{enumerate}
\end{example}

\begin{example}
$F(x,y,z)=x^4+y^2+yz^2-z^2$, $G(x,y,z)=xy$. $\nabla
F(0,0,0)=(0,0,0)$, $\nabla G(0,0,0)=(0,0,0)$. It is easy to know
that the origin $P=(0,0,0)$ is a 2-fold point of both surfaces
$F(x,y,z)=0$ and $G(x,y,z)=0$, then $P$ is a 4-fold point of the
curve since the tangent planes of the two surfaces are different (as
shown in Fig.~\ref{fig:S-curve-4}).

By Algorithm \ref{Alg:ql3} and a parametric cubic curve
(\ref{eq:PC3}), we obain
\begin{eqnarray*}
\left\{\begin{array}{ll}F(x,y,z)=C_{1}t+C_{2}t^2+C_{3}t^3+\cdots,\\
G(x,y,z)=D_{1}t+D_{2}t^2+D_{3}t^3+\cdots,\end{array}\right.
\end{eqnarray*}
where, $C_1=0,\ D_1=0,\ C_2=b_1^2-c_1^2,\ D_2=a_1b_1$.

By $C_2=D_2=0$, we obtain four tangents $(0,1,1),\ (0,1,-1),\
(1,0,0)$ and $(1,0,0)$.
\begin{enumerate}[1)]
\item If $a_1=0,\ c_1=b_1$,
the tangent to the curve at $P$ is $(0,1,1)$, and
$C_3=b_1^3+b_1b_2-b_1c_2$, $D_3=a_2b_1/2$.

Let $C_3=D_3=0$, we obtain $a_2=0,\ c_2=b_1^2+b_2$. By
Eq.~(\ref{eq:ql2}), the curvature of the curve at $P$ is
$\frac{1}{2\sqrt{2}}$. Let $C_4=D_4=0$, then by Eq.~(\ref{eq:nl}),
the torsion of the curve at $P$ is $0$.

\item If $a_1=0,\ c_1=-b_1$, we obtain the curvature and torsion of the curve at $P$
along the tangent $(0,1,-1)$ are $\frac{1}{2\sqrt{2}}$ and 0.

\item If $b_1=c_1=0$, the tangent to the curve at $P$ is $(1,0,0)$ (a double root), and $C_3=0$,
$D_3=a_1b_2/2$, $C_4=a_1^4+b_2^2/4-c_2^2/4$,
$D_4=a_2b_2/4+a_1b_3/6$. The curvature cannot be determined only by
$C_3=D_3=0$. Note that $P$ is a 4-fold point of the curve, it is
necessary that $C_3=D_3$, $C_4=D_4$. Then we have $b_2=0,\
b_3=\frac{3(4a_1^4-c_2^2)}{2a_1}$, and $C_3=D_3=0$,
$C_4=D_4=a_1^4-c_2^2/4$. By $C_4=D_4=0$ and Eq.~(\ref{eq:ql2}), the
curvature of the curve at $P$ is $2$ (of multiplicity 2). Besides,
we obtain $b_3=0$, then $y=0$. It means that the curve lies on a
plane, hence the torsion is $0$.
\end{enumerate}
\end{example}

\end{CJK*}

\end{document}